\documentclass[11pt,leqno]{article}
\usepackage{amssymb,amsmath,amscd,amsthm,amsfonts,graphicx}
\usepackage[dvips]{pstcol}
\usepackage[all]{xy}
\usepackage{epsfig}

\newtheorem{theorem}{Theorem}[section]
\newtheorem{lemma}[theorem]{Lemma}
\newtheorem{proposition}[theorem]{Proposition}
\newtheorem{corollary}[theorem]{Corollary}
\newtheorem{remark}[theorem]{Remark}
\newtheorem{definition}[theorem]{Definition}

\newtheorem{example}[theorem]{Example}

\def\semi{\hbox{ $\times $ \kern-.972em \raise.12719em\hbox{ $_{^|}$}  }}

\def\be{\begin{enumerate}}
\def\ee{\end{enumerate}}
\def\bi{\begin{itemize}}
\def\ei{\end{itemize}}

\title{Conjugacy in Garside Groups II: \\ Structure of the Ultra Summit Set}
\author{Joan S. Birman\footnote{Partially
supported by the U.S.National Science
Foundation, under Grant DMS-0405586.}, Volker Gebhardt,
Juan Gonz\'alez-Meneses\footnote{Partially supported by MTM2004-07203-C02-01 and FEDER.}.}
\date{June 26, 2006.}

\begin{document}
\maketitle

  \begin{abstract}
  \noindent This paper is the second in a series (the others are \cite{BGGM-I} and \cite{BGGM-III}) in which the authors study the conjugacy decision problem (CDP) and the conjugacy search problem (CSP) in Garside groups.

  The ultra summit set $USS(X)$ of an element $X$ in a Garside group $G$ is a finite set of elements in $G$, introduced in~\cite{gebhardt}, which is a complete invariant of the conjugacy class of $X$ in $G$. A fundamental question, if one wishes to find bounds on the size of $USS(X)$, is to understand its structure. In this paper we introduce two new operations on elements $Y\in USS(X)$, called `partial cycling' and `partial twisted decycling' (Definition \ref{D:partial cycling}),  and prove that if $Y, Z\in USS(X)$, then $Y$ and $Z$  are related by sequences of partial cyclings and partial twisted decyclings.  These operations are a concrete way to understand the minimal simple elements which result from the convexity theorem in \cite{gebhardt}.   Using partial cycling and partial twisted decycling, we investigate the structure of a directed graph $\Gamma_X$ which is determined by $USS(X)$, and show that $\Gamma_X$ can be decomposed into `black' and `grey'  subgraphs.  There are applications which relate to the program outlined in \cite{BGGM-I} for finding a polynomial solution to the CDP/CSP in the case of braids. A different application is to give a new algorithm for solving the CDP/CSP in Garside groups which is faster than all other known algorithms, even though its theoretical complexity is the same as that of the algorithm in \cite{gebhardt}.  There are also applications to the theory of reductive groups.
\end{abstract}

\setlength{\parskip}{0 pt}
\tableofcontents
\setlength{\parskip}{1ex plus 0.5ex minus 0.2ex}

\section{Introduction}
\label{section:Introduction}

Garside groups were introduced by P. Dehornoy and L. Paris in~\cite{D-P}.  Roughly speaking, they are groups that satisfy the same kind of
properties that were shown to hold in Artin's braid groups $\{B_n, \ n=1,2,\dots\}$ in the foundational paper of F. Garside~\cite{garside}, and soon afterwards in spherical Artin-Tits groups by
Brieskorn and Saito~\cite{B-S}.  More
precisely, they are groups admitting a lattice structure, invariant under
left multiplication, and satisfying some additional conditions.
The focus of this paper and of its companions \cite{BGGM-I} and  \cite{BGGM-III}  is on gaining an improved understanding of the known algorithmic solutions in Garside groups of the {\it conjugacy decision problem} (CDP),
that is, to decide, for arbitrary $X,Y\in G$ whether there exists $Z$ such that $Y = Z^{-1}XZ$; and also the {\it conjugacy search problem} (CSP), that is,  to find $Z$ as above when $X$ and $Y$ are known to be conjugate.

In the manuscript \cite{gebhardt} the second author of this paper introduced a complete invariant of the conjugacy class of an arbitrary element $X$ in an arbitrary Garside group $G$, which is a certain finite set of elements, called the `ultra summit set' $USS(X)$ (Definition \ref{D:USS}.)  For arbitrary  $X,Y\in G$ either $USS(X) = USS(Y)$ (in which case $X$ is conjugate to $Y$) or $USS(X)\cap USS(Y) = \emptyset$.  From this it follows that in order to determine whether $X$ and $Y$ are conjugate, one must compute the entire set $USS(X)$ and a single element in $USS(Y)$.   That leads one immediately to ask about bounds on the number of elements in $USS(X)$. This paper is a first step toward answering that question. The focus of this paper is on the structure of $USS(X)$, a matter which we isolated as being a key issue as we worked toward a polynomial algorithm for the CDP/CSP in $B_n$.   See $\S$1.4 of \cite{BGGM-I} for an outline of our approach to that problem, and a discussion of how the work in this paper fits into that program and the open problems that remain.

Here is a brief overview of the contents of this paper.  In $\S$\ref{SS:basic facts} we give a quick review of the necessary background, with all essential defintions.   In $\S$\ref{S:the structure of the USS} we will describe new properties of the structure of ultra
summit sets, showing that in order to compute them one only needs to
apply two special kinds of conjugation, which we call `partial
cycling' and `partial twisted decycling' (Definition \ref{D:partial cycling}).   Using it, we discover new aspects
 of the combinatorics that are both closely related to things we already understand well, but at the same time present new challenges.  We will study the graph $\Gamma_X$ which describes $USS(X)$ of an element $X\in G$, a Garside group, and prove that $\Gamma_X$ splits into two subgraphs, which we call `black' and `grey' .  The similarities and differences between the graphs $\Gamma_X$ and $\Gamma_{X^{-1}}$ are particularly interesting.  In $\S$2 we also give several key examples which illustrate our new discoveries.

In $\S$\ref{S:rigid elements} and $\S$\ref{S:periodic elements} we examine the effects of partial cycling on two special kinds of elements in a Garside group $G$: rigid elements and periodic elements.  In this regard we note that, while the work in this paper is applicable to all Garside groups and no assumptions are made that restrict attention to the braid groups $B_n$, the reader who is familiar with the Thurston-Nielsen trichotomy (see $\S$1.3 of \cite{BGGM-I}) in braid groups and with the main result in \cite{BGGM-I} will recognize that rigid elements and periodic elements have particular importance in the case of braids.  This is relevant for the program  that we outlined in $\S$1.4 of \cite{BGGM-I} for studying the CDP/CSP in the braid groups.

The final section, $\S$\ref{S:applications} is devoted to applications.  The first application, in
$\S$\ref{SS:complexity of the CDP/CSP for PA and periodic braids}, is to show how our new knowledge of the structure of $USS(X)$ in $\S$\ref{S:rigid elements} and $\S$\ref{S:periodic elements} advances our efforts to find a polynomial solution to the CDP/CSP in the special case of the braid groups. Unfortunately, we were unable to find a bound  for the number of elements in $USS(X)$ in this paper,  even when we restrict to  rigid elements.  Finding such a bound is open question 2 of $\S$1.4 of \cite{BGGM-I}, and remains for future work.  Examples are given in $\S$\ref{SS:complexity of the CDP/CSP for PA and periodic braids} that illustrate the remaining difficulties.  Such examples would have been impossible to understand without the work in this paper.   We then discuss the periodic case, briefly. A polynomial solution to the CDP/CSP for periodic braids will be given in the next paper in this series, \cite{BGGM-III}.

The second application, in $\S$\ref{SS:a new solution to CDP/CSP in Garside groups} is to develop a new solution to the CDP/CSP in arbitrary Garside groups. This new algorithm is faster than the one in~\cite{gebhardt}, although its theoretical complexity is essentially the same.   The third application, in $\S$\ref{SS:partial cycling subsumes decyling} is to show that, in a sense, partial cycling subsumes cycling. This is probably more of theoretical interest than of interest for complexity.  The final application, in $\S$\ref{SS:an application to the theory of reductive groups}, is a brief note to mention that, as was communicated to us by Jean Michel,  our work on the structure of the $USS(X)$ of a periodic element has implications for the theory of reductive groups.   We refer to \cite{DM} for details.

\subsection{Definitions and basic facts}
\label{SS:basic facts}

Our review of what we need from the literature will be brief.  We strongly recommend that non-experts consult $\S$1  of \cite{BGGM-I} for a more leisurely introduction to Garside groups and a review of what is known about the CDP/CSP with lots of examples.

\begin{definition}
\label{D;Garside group} {\rm A group $G$ is said to be a {\it Garside group} (of finite type) if it satisfies the following three properties:
\be
\item [(A)] $G$ admits a lattice order $(G,\preccurlyeq,\vee,\wedge)$, invariant under
left-multiplication.

This means in particular that, with respect to the order $\preccurlyeq$, every pair of elements $s,t\in G$  admits a unique lcm
$s\vee t$ and a unique gcd $s\wedge t$.
Let $P=\{p\in G;\; 1\preccurlyeq p\}$.  Then $P\subset G$ is a monoid.   We call the elements of $P$ the {\it positive} elements of $G$.

\item [(B)] There exists an element $\Delta \in P$, called the {\it Garside element}, satisfying two properties: (a) The interval $[1,\Delta]=\{s\in G;\; 1\preccurlyeq s \preccurlyeq \Delta\}$ is finite and
generates $G$. (b)$\Delta^{-1} P \Delta = P$.

\item [(C)]  The monoid $P$ is {\it atomic}, that is $\forall X\in P$ there exists an upper bound on the length of a chain $1 \prec X_1 \prec \cdots \prec X_n = X$.  The {\it atoms} are those elements $a\in P$ for which $u\in P, \ \ 1\neq u\preccurlyeq a \Longrightarrow u=a$.  The atoms are in $[1,\Delta]\subset P$ and they generate $G$.
\ee}
\end{definition}

Many examples of Garside groups can be found in \cite{BGGM-I}.

\begin{definition}
\label{D:simple prefix}
{\rm  The elements in $[1,\Delta]$ are called the {\it simple elements} of $G$. A {\it prefix} of $x\in G$
is an element $q\in G$ such that $q^{-1}x\in P$.
A {\it simple prefix} $q$ of $x$ is a prefix which is simple.}
\end{definition}

\begin{definition}
\label{D:left normal form}
{\rm Given $X\in G$, we say that a decomposition of $X$ is in {\it left normal form}:}
\begin{equation}
\label{E:left normal form}
X=\Delta^p
x_1\cdots x_r \ \  (r\geq 0)
\end{equation}
{\rm if  $p\in \mathbb Z$ is the maximal  power such that $\Delta^p\preccurlyeq X$, and  $x_i\not= 1$ is the maximal simple prefix of $x_ix_{i+1}\cdots x_r$ for every $i=1,\dots,r$.   A pair of
simple elements $a,b\in [1,\Delta]$ is {\it left weighted}  if the product $ab$ is
in left normal form as written.  If $X=\Delta^p x_1\cdots x_r$ is in left normal form, the {\it
infimum, supremum} and {\it canonical length} of $X$, are defined by
$\inf(X)=p$, $\sup(X)=p+r$ and $\ell(X)=r$, respectively. }
\end{definition}

\begin{definition}  {\rm The symbol $\tau$ denotes the inner automorphism of $G$ that is defined by $\tau(x)=\Delta^{-1} x \Delta$. It is called the {\it twisting automorphism}.}
\end{definition}

\begin{definition}  {\rm  If $u,v\in G$ and $uv=\Delta$, then $v=\partial(u)$.   It is a basic property of the lattice structure on $G$ that  the map $u\to\partial(u)$ is a permutation of the finite set of simple elements.   Therefore if $u$ is simple and $v=\partial(u)$, it makes sense to say that $u=\partial^{-1}(v)$. }
\end{definition}

Let us recall (see \cite{el-m}) that if $X$ has left normal form as in (\ref{E:left normal form}), then the left normal form of $X^{-1}$ is determined by that of $X$:
\begin{equation}
\label{E:left normal form of inverse}
X^{-1} = \Delta^{-p-r}x_r^\prime\cdots x_1^\prime, \ {\rm where} \ \ x_i^\prime = \tau^{-p-i}(\partial(x_i)).
\end{equation}
This interplay between $X$ and $X^{-1}$ will play an important role in our work.

There is a special type of conjugation that has been recognized, ever since the
work of ElRifai and Morton in \cite{el-m},  to be of fundamental importance in
the combinatorics that underlie all of our work:

\begin{definition} \label{D:cycling,decycling}
{\rm If $\ell(X) >0$ and has left normal form as in (\ref{E:left normal form}), then {\it cycling}  $\mathbf c(X)$ and {\it decycling} $\mathbf d(X)$ are the special conjugates of $X$ that are defined by:}
\begin{equation}
\label{E:cyclings and decyclings}
\mathbf c(X)  = \Delta^p x_2\cdots x_r  \:\tau^{-p}(x_1) \ \ \ \ \ \ \ \ \ \ \  {\rm and}  \ \ \ \ \ \ \ \ \ \ \ \  \mathbf d(X)  = x_r  \Delta^p  x_1\cdots x_{r-1}
\end{equation}
{\rm The {\it twisted decycling} of $X$ is defined by $\tau(\mathbf d(X))$.}
\end{definition}

Twisted decycling will be important later.  To see how cycling and decycling enter into the picture,  recall that the {\it super summit set} $SSS(X)$ of $X\in G$ is the set of all elements $Y $ in the conjugacy class of $X$ with the property that if $Y = \Delta^q y_1\cdots y_s$ in left normal form, then $q = \inf(Y)$ is maximal and $q+s = \sup(Y)$ is minimal.  It is known that $SSS(X)$ is always non-empty and finite.
See Corollary 1.11 of \cite{BGGM-I}  and the discussion that follows it for a description of how cycling and decycling enable one, for any $X\in G$, to compute at least one element  $\widetilde{X} \in SSS(X)$. However, we will not be working with $SSS(X)$ in this paper.  Instead, we consider a subset of $SSS(X)$ which was introduced by the second author of this paper in~\cite{gebhardt}.
\begin{definition}
\label{D:USS}
{\rm The {\it ultra summit set} $USS(X)$ of an element $X\in G$ is  the subset of elements $Y\in SSS(X)$ such that $\mathbf
c^m(Y)=Y$, for some $m > 0$. }
\end{definition}

Note that $USS(X)$  is always non-empty, because if $Y\in SSS(X)$, then
${\bf c}(Y)$ gives us another element of $SSS(X)$, and since $SSS(X)$ is
finite, it follows that that after some number of iterations of cycling we
will find integers $m_1, m_2$, with $m_1<m_2,$ such that ${\bf c}^{m_1}(Y) =
{\bf c}^{m_2}(Y)$.   But then ${\bf c}^{m_1}(Y)$ is in $USS(X)$.

\begin{definition}
\label{D:minimal simple element}
{\rm Given $X\in G$ and $Y\in USS(X)$, we say that a simple element $s\not= 1$
is a {\it minimal simple element} for $Y$ with respect to $USS(X)$
if $Y^s=s^{-1}Ys \in USS(X)$, and no proper prefix of $s$ satisfies
this property. }
\end{definition}

In~\cite{gebhardt} it is shown how to compute the minimal simple
elements corresponding to a given $Y\in USS(X)$, and this allows one to
compute the whole $USS(X)$, starting with a single element $\widetilde X\in USS(X)$.  To state the theorem that we will use from~\cite{gebhardt}, we need one more concept:

\begin{definition}
\label{D:Gamma_X}
{\rm Given $X\in G$, there is a directed graph, $\Gamma_X$, which describes the entire set $USS(X)$, is defined by the
following data:}
\be
 \item  {\rm Each vertex represents an element $Y\in USS(X)$.}
 \item {\rm For every $Y\in USS(X)$ and every minimal simple element $s$
 for $Y$ with respect to $USS(X)$, there is an arrow labeled by $s$
 going from $Y$ to $Y^s$.}
\ee
\end{definition}

The algorithm in~\cite{gebhardt} computes the graph $\Gamma_X$ using
the following `convexity theorem', which is analogous to a related convexity theorem proved for super
summit sets in \cite{F-GM}.  It will be seen to be of fundamental importance in our work, and will be given concrete meaning in many examples.

\begin{theorem}\label{T:convexity_USS} {\rm \cite{gebhardt}}
Let $X\in G$ and $Y\in USS(X)$.
{\rm (i)}~If $s,t\in G$ are such that $Y^s\in USS(X)$ and
$Y^t\in USS(X)$, then $Y^{s\wedge t}\in USS(X)$.
{\rm (ii)}~For every $u\in P$ there is a unique element $c_Y(u)$  which is minimal with respect to:
$ u\preccurlyeq c_Y(u)$ and  $Y^{c_Y(u)}\in USS(X).$
{\rm (iii)}~The graph $\Gamma_X$ that is described in Definition~\ref{D:Gamma_X} is connected.
\end{theorem}

Notice that the set of minimal simple elements for $Y$ is contained
in $c_Y(A) = \{c_Y(a)\; | \; a \mbox{ is an atom}\}$, hence the
number of minimal simple elements for $Y$ is bounded by the number
of atoms.  In the case of $B_n$, with its classical Garside structure, there are $n-1$ atoms, so in
$\Gamma_X$ there are at most $n-1$ arrows starting at a given vertex.

We remark that one obtains an element in $USS(X)$ by
iterated application of cycling to an element in $SSS(X)$. The number of
times one needs to apply cycling, in order to go from an element in
$SSS(X)$ to an element in $USS(X)$ is not known in general.  That is open question 3 in $\S$1.4 of \cite{BGGM-I}.
Nevertheless, the theoretical complexity of the algorithm
in~\cite{gebhardt} is not worse than the theoretical complexity of the algorithm
in~\cite{F-GM}, which is based upon $SSS(X)$, and the performance of the algorithm is substantially better in practice.  Later in this paper, in
$\S$\ref{SS:a new solution to CDP/CSP in Garside groups}, we will obtain yet another algorithm which is faster than the one in \cite{gebhardt}, although its theoretical complexity is unchanged.

The conjugating elements involved in
cycling or decycling will play a crucial role in our work (as they did in \cite{BGGM-I}), so to avoid repeated references to the automorphism $\tau$ we introduce convenient shorthand:
\begin{definition} \label{D:iota_varphi}
{\rm Given $X\in G$ whose left normal form is $X=\Delta^p x_1\cdots
x_r$ ($r>0$), we define the {\it initial factor} of $X$ as
$\iota(X)=\tau^{-p}(x_1)$, and the {\it final factor} of $X$ as
$\varphi(X)=x_r$. If $r=0$ we define $\iota(\Delta^p)=1$ and
$\varphi(\Delta^p)=\Delta$.
Equivalently $ \iota(X) =
X\Delta^{-p}\wedge\Delta$ and $\varphi(X) = (\Delta^{p+r-1}\wedge X)^{-1}X.$}
\end{definition}

Rigid elements in a Garside group were introduced in \cite{BGGM-I}:
\begin{definition}
\label{D:rigid}
{\rm Let $X = \Delta^p x_1\cdots x_r$ be in left normal form, with
$r = \ell(X) >0$. Then $X$ is {\it rigid} if the element $\Delta^p\: x_1\cdots
x_r  \:\tau^{-p}(x_1)\, $ is in left normal form as written.}
\end{definition}

Notice that if $X$ is rigid, then the cycling of $X$, that is,
$\mathbf c(X) = \Delta^p\: x_2\cdots x_r
 \:\tau^{-p}(x_1)\, $ is in left normal form as written. Actually,
this latter property is equivalent to $X$ being rigid if $r>1$. But
we prefer the definition above, otherwise every element of canonical
length 1 would be rigid.

\begin{proposition} {\rm \cite{BGGM-I}}
\label{P;equivilent defs of rigidity} The following are equivalent
characterizations of rigid elements. {\rm (i)} $X$ is rigid. ~{\rm
(ii)} $\ell(X)>0$ and $\varphi(X)\iota(X)$ is left weighted as
written. In particular normal forms are preserved under cycling.
~{\rm (iii)} $\ell(X)>0$ and  $\iota(X)\wedge \iota(X^{-1})=1$.
\end{proposition}

The third condition implies that $X$ is rigid if and only if $X^{-1}$ is rigid. It is important to notice that the assumption $\ell(X)>0$ implies that $\Delta$ and its powers are {\it not}  rigid.

\begin{remark} \label{R:T-N} {\rm We need one more set of ideas from the existing literature before we can begin our new work.   In $\S$\ref{SS:complexity of the CDP/CSP for PA and periodic braids} we will specialize our work about Garside groups to the braid groups $B_n$.  In addition, all of our examples will be taken from braid groups.
It therefore is important to keep in mind that the braid group $B_n$ acts on the $n$-times punctured disc $D_n$, and that this action determines a faithful representation of $B_n$ as the {\it mapping class group}  of $D_n$, that is, the group $\pi_0({\rm Diff}^+(D_n))$. The reason this is important is that, as a consequence, the well-known {\it Thurston-Nielsen Trichotomy}, introduced by Thurston in \cite{thurston} applies: braids come in three flavors: pseudo-Anosov (or PA), periodic and reducible.  We refer the reader to \cite{thurston} for the foundational paper in this area and to \cite{BGGM-I} for a very brief description of the  trichotomy.   We will not discuss it here because we use it in a very peripheral way: to identify our examples as being PA or periodic.  We note, however, that in~\cite{BGGM-I} the authors of this paper proved the important fact that every PA braid in its $USS$ has a `small' rigid power, so that all results about rigid elements in Garside groups apply ipso-facto to powers of PA braids.  A major part of the work in this paper, notably the work in $\S$\ref{S:rigid elements} concerns rigid elements in Garside groups.  As a result, much of the work in this paper has special relevance in the key case of PA braids. }
\end{remark}

{\bf Acknowledgements:}  The authors thank Jean Michel, who showed us the importance of partial cyclings in the study of reductive groups, and encouraged us to show the results in $\S$\ref{S:periodic elements}.

J. Birman thanks the Project MTM2004-07203-C02-01 of the Spanish
Ministerio de Ciencia y Tecnolog\'\i a for hosting her visit to
Seville in November 2004, so that she and J. Gonzalez-Meneses could
work together on this project.

J. Gonz\'alez-Meneses thanks the project MTM2004-07203-C02-01 and
the Columbia University Department of Mathematics for hosting his
two visits to New York, in July 2004 and March-April 2006.

V. Gebhardt thanks the department of Algebra of the University of
Seville, and the Junta de Andaluc\'\i a, for funding his visit to
Seville in January 2006.

\section{The structure of the ultra summit set}
\label{S:the structure of the USS}
In this section we begin the new work in this paper. In $\S$\ref{SS:description of the minimal simple elements} we begin to investigate the minimal simple elements whose existence was established in statement (ii) of Theorem~\ref{T:convexity_USS}. In particular,  in
Theorem~\ref{T:minimal_simple_elements} we give new meaning to the somewhat elusive element $c_Y(u)$.  Using our description, we introduce partial cycling  and partial twisted decycling in Definition \ref{D:partial cycling}.   The main result in $\S$\ref{SS:description of the minimal simple elements} is Corollary~\ref{C:partial cycling and USS(X)}.  As a result of this Corollary, we will prove that $\Gamma_X, \ X\in G$ splits into two (in general not connected) subgraphs which we call `grey' and `black'.   The properties of the grey and black subgraphs are studied in $\S$\ref{SS:black_grey} and $\S$\ref{SS:Intersection_of_components}.

\subsection{Description of the minimal simple elements}
\label{SS:description of the minimal simple elements}
As was proved in  \cite{BGGM-I}, the initial and final factors $\iota(X), \varphi(X)$ and $\iota(X^{-1}), \
\varphi(X^{-1})$ are closely related:

\begin{lemma}\label{L:initial & final factors of inverses} {\rm \cite{BGGM-I}}
For every $X\in G$ one has $\iota(X^{-1})=\partial(\varphi(X))$ and
$\varphi(X^{-1})=\partial^{-1}(\iota(X))$.
\end{lemma}

We saw in Theorem~\ref{T:convexity_USS} that, given one element
$X\in USS(X)$, one can obtain any other element in $USS(X)$ just
through conjugations by minimal simple elements. In this section we will describe the
minimal simple elements for $X\in USS(X)$, and show how they determine the graph $\Gamma_X$.  First, we need some results.

\begin{lemma}\label{L:conjugating elements}
Given $X\in G$ with $\ell(X)>0$, one has:
$$
X^{\iota(X)}= \mathbf c(X), \hspace{1cm} X^{\varphi(X)^{-1}}= \mathbf d(X), \hspace{1cm} X^{\iota(X^{-1})}= X^{\partial(\varphi(X))}= \tau(\mathbf d(X)).
$$
\end{lemma}

\begin{proof}
The first two claims follow from the definitions. The third one is shown as follows: $X^{\partial(\varphi(X))} = X^{\varphi(X)^{-1} \Delta} =(\mathbf d(X))^\Delta = \tau (\mathbf d(X))$.
\end{proof}

\begin{remark}
\label{twisted decycling}
{\rm By the above result, the cycling and the twisted decycling of $X$ are conjugates of $X$ by simple elements. Moreover, the twisted decycling of $X$ corresponds to a cycling of $X^{-1}$, since it is a conjugation by $\iota(X^{-1})$. In the same way, the cycling of $X$ corresponds to a twisted decycling of $X^{-1}$.}
\end{remark}

\begin{lemma}\label{L:tau_c_d}
The twisting automorphism $\tau$ commutes with both cycling and decycling. Moreover, the USS of an element is closed under twisting, cycling and decycling.
\end{lemma}

\begin{proof}
It is clear from the definitions that
$\iota(\tau(X))=\tau(\iota(X))$ and
$\varphi(\tau(X))=\tau(\varphi(X))$. In particular, $\mathbf
c(\tau(X))= \tau(X)^{\tau(\iota(X))}= \tau(X^{\iota(X)})=
\tau(\mathbf c(X))$, and $\mathbf d(\tau(X))=
\tau(X)^{\tau(\varphi(X)^{-1})}= \tau(X^{\varphi(X)^{-1}})=
\tau(\mathbf d(X))$, that is, $\tau$ commutes with both cycling and
decycling.

From now on let $X\in USS(X)$.  We have to show that $\tau(X)$,
$\mathbf c(X)$ and $\mathbf d(X)$ belong to $USS(X)$. Since $X\in
USS(X)\subseteq SSS(X)$, and $\tau(X)$ has the same canonical length
as $X$, it follows that $\tau(X)\in SSS(X)$. Moreover, one has
$\mathbf c^m(X)=X$ for some $m$.  Then, $\mathbf
c^m(\tau(X))=\tau(\mathbf c^m(X))=\tau(X)$ and therefore $\tau(X)\in
USS(X)$.

As $\ell(\mathbf c(X))\le\ell(X)$ \cite{el-m} and $X$ has minimal canonical length
among its conjugates, we have $\ell(\mathbf c(X))=\ell(X)$, that is,
$\mathbf c(X)\in SSS(X)$.  It is then clear that $\mathbf c(X)\in
USS(X)$, since $\mathbf c^{m}(\mathbf c(X))=\mathbf c(\mathbf c^{m}(X))=\mathbf c(X)$.

Finally, notice that $X^X=X\in USS(X)$ and that
$X^{\Delta^{p+r-1}} = \tau^{p+r-1}(X) \in USS(X)$. Then, by
Theorem~\ref{T:convexity_USS}, $\mathbf d(X) = X^{(\Delta^p
x_1\cdots x_{r-1})} = X^{X\wedge \Delta^{p+r-1}} \in USS(X)$.
\end{proof}

\begin{theorem}\label{T:minimal_simple_elements}
Let $X\in USS(X)$ and let $s$ be a minimal simple element for
$X$. Then one and only one of the following conditions holds:
\begin{enumerate}
\item $\varphi(X)s$ is a simple element.
\item $\varphi(X)s$ is left weighted as written.
\end{enumerate}
\end{theorem}

\begin{proof}
Recall from Lemma~\ref{L:conjugating elements} that $X^{\partial(\varphi(X))} = X^{\varphi(X)^{-1} \Delta} =
(\mathbf d(X))^\Delta = \tau (\mathbf d(X))$. Hence, by
Lemma~\ref{L:tau_c_d}, $X^{\partial(\varphi(X))}\in USS(X)$. On the
other hand, since $s$ is a minimal simple element for $X$, one has
$X^s \in USS(X)$. Therefore, by Theorem~\ref{T:convexity_USS} it
follows that $X^{\partial(\varphi(X)) \wedge s}\in USS(X)$.

Consider $t= \partial(\varphi(X))\wedge s$.  By definition,
$t\preccurlyeq s$, and we just saw that $X^t\in USS(X)$. Since $s$ is a
minimal simple element for $X$, this implies that either $t=s$ or
$t=1$.

Notice that the first factor in the left normal form of $\varphi(X) s$
(which is possibly $\Delta$) is equal to
$\varphi(X) s \wedge \Delta
 = \varphi(X)s \wedge \varphi(X) \partial(\varphi(X))
 = \varphi(X) (s \wedge \partial(\varphi(X)))
 = \varphi(X)t$.
Hence it is equal to either $\varphi(X)s$ or $\varphi(X)$. The first case
implies that $\varphi(X)s$ is simple, while the second one means
that $\varphi(X)s$ is left weighted.
\end{proof}

\begin{proposition}
\label{P:prefix} Let $X\in USS(X)$ with $\ell(X)>0$ and let $s$ be a
minimal simple element for $X$ such that $\varphi(X)s$ is left
weighted.  Then, $s\preccurlyeq \iota(X)$.
\end{proposition}
\begin{proof}
Let $\Delta^p x_1\cdots x_r$ be the left normal form of $X$. If
$\varphi(X)s=x_r s$ is left weighted, then $\Delta^p x_1\cdots x_r
s$ is the left normal form of $X s$. However, we know that $s^{-1}
Xs = \Delta^p \tau^p(s)^{-1} x_1\cdots x_r s \in USS(X)$, hence
$\tau^p(s)^{-1} x_1\cdots x_r s\in P$, which is equivalent to
$\tau^p(s)\preccurlyeq x_1\cdots x_r s$. Since $\tau^p(s)$ is simple,
this in turn is equivalent to $\tau^p(s)\preccurlyeq x_1\cdots x_r
s\wedge\Delta$. Finally, since $x_1\cdots x_r s$ is in left normal
form with $r>0$ , this means that $\tau^p(s)\preccurlyeq x_1$, thus
$s\preccurlyeq \tau^{-p}(x_1) = \iota(X)$, as we wanted to show.
\end{proof}

\begin{corollary}\label{C:prefix}
Let $X\in USS(X)$ with $\ell(X)>0$ and let $s$ be a minimal simple
element for $X$. Then $s$ is a prefix of either $\iota(X)$ or
$\iota(X^{-1})$, or both.
\end{corollary}

\begin{proof}
We have seen in Lemma~\ref{L:initial & final factors of inverses} that
$\iota(X^{-1})=\partial(\varphi(X))$ and we know by
Theorem~\ref{T:minimal_simple_elements} that $\varphi(X)s$ is either
simple or left weighted. In the first case, $s\preccurlyeq
\partial(\varphi(X)) = \iota(X^{-1})$, whereas in the second case
$s\preccurlyeq \iota(X)$, by Proposition~\ref{P:prefix}.

Notice that we could have at the same time $s\preccurlyeq \iota(X)$ and
$s\preccurlyeq \iota(X^{-1})$, but this is only possible in the first
case, where $\varphi(X)s$ is simple.
\end{proof}

We remark that the statements of Proposition~\ref{P:prefix} and
Corollary~\ref{C:prefix} are wrong for elements $X$ with
$\ell(X)=0$, as in this case $\iota(X)=\iota(X^{-1})=1$.
\medskip

Recall that every two elements in $USS(X)$ can be joined by a
sequence of conjugations by minimal simple elements. We have just
seen that the minimal simple elements for each $Y\in USS(X)$ with
$\ell(Y)>0$ are prefixes of either $\iota(Y)$ or
$\iota(Y^{-1})=\partial(\varphi(Y))$. Let us give a name to these
special kinds of conjugations.

\begin{definition}
\label{D:partial cycling}
{\rm Let $X\in G$.  A {\it partial cycling of  $X$} is a conjugation of
$X$ by a prefix of $\iota(X)$. A {\it partial twisted decycling of $X$} is a conjugation of $X$ by a prefix of $\iota(X^{-1})$, or equivalently (using Lemma~\ref{L:initial & final factors of inverses}), by a prefix of $\partial(\varphi(X))$.}
\end{definition}

\begin{corollary}
\label{C:partial cycling and USS(X)}
Given $X, Y\in USS(X)$, there exists
a sequence of partial cyclings and partial twisted decyclings  joining $X$ to $Y$.
\end{corollary}

\begin{proof}
This is a direct consequence of Theorem~\ref{T:convexity_USS} and
Corollary~\ref{C:prefix}.
\end{proof}

We end this section with several examples which will be used now and later in this paper to illustrate various points.   All of our examples will be taken from braid groups, with the classical Garside structure that was discovered in \cite{garside}, and with the elementary braids $\sigma_1,\dots,\sigma_{n-1}$ as atoms.
By Corollary~\ref{C:prefix}, there are two kinds
of minimal simple elements, hence there are two kinds of arrows in
$\Gamma_X$. We say that an arrow $s$ starting at a vertex
$Y=\Delta^p y_1\cdots y_r$ is {\it black} if $s$ is a prefix of
$\iota(Y)$, and {\it grey} if $s$
is a prefix of $\iota(Y^{-1})$ or, equivalently, if $y_r s$ is a
simple element. In other words, an arrow starting at $Y$ is black if
it corresponds to a partial cycling of $Y$, and it is grey if it corresponds to a partial twisted decycling of $Y$.

\begin{figure}[htpb]
\centerline{\includegraphics[scale=1.0]{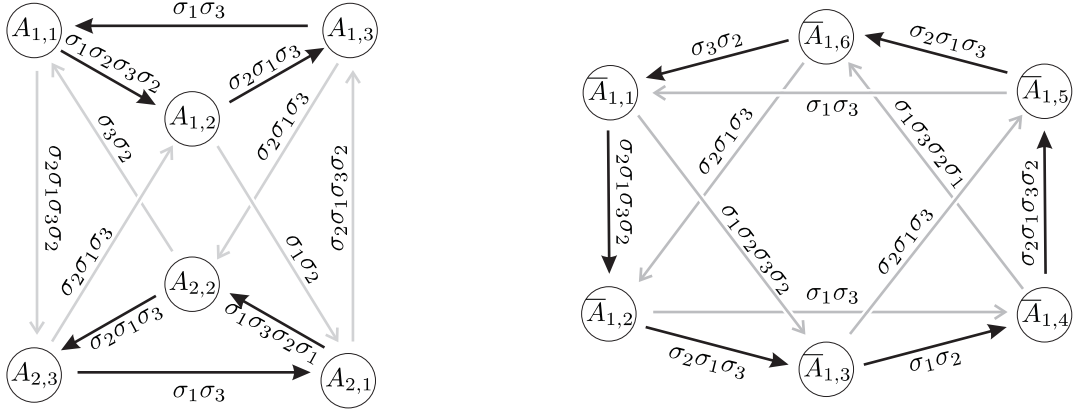}} \caption{The graphs
$\Gamma_A$ and $\Gamma_{\overline A} = \Gamma_{A^{-1}}$.}
\label{F:A}
\end{figure}

\begin{example}
\label{E:A} {\rm Our first example is the 4-braid, $A =  \sigma_1
\sigma_2 \sigma_3 \sigma_2 \sigma_2 \sigma_1 \sigma_3 \sigma_1
\sigma_3$.   The reader is referred to Remark~\ref{R:T-N} and thence
to $\S$1.3 of \cite{BGGM-I} and \cite{thurston} for a brief
description of the so-called {\it Thurston-Nielsen trichotomy}  in
mapping class groups.  The braid $A$ is PA  and rigid (see
Remark~\ref{R:T-N} and Definition~\ref{D:rigid}).  Its first cycling
orbit $A_1$ has length 3, and by Lemma~\ref{L:tau_c_d}, since
$\tau(A_1)\not= A_1$ there is a second orbit  $A_2 = \tau(A_1)$. A
computation shows that $USS(A)$ has exactly these two cycling
orbits, with 3 elements each, namely:}
\begin{eqnarray}
\nonumber  A_1 = \{A_{1,1} &=& \sigma_1 \sigma_2 \sigma_3 \sigma_2 \cdot \sigma_2 \sigma_1 \sigma_3 \cdot \sigma_1 \sigma_3,  \\
\nonumber  A_{1,2} &=& \sigma_2 \sigma_1 \sigma_3 \cdot \sigma_1 \sigma_3 \cdot \sigma_1 \sigma_2 \sigma_3 \sigma_2, \\
\nonumber  A_{1,3} &=&  \sigma_1 \sigma_3\cdot \sigma_1 \sigma_2
\sigma_3 \sigma_2 \cdot \sigma_2 \sigma_1 \sigma_3.\}
\end{eqnarray}
\begin{eqnarray}
 \nonumber  A_2 = \{A_{2,1} &=&  \sigma_1 \sigma_3 \sigma_2 \sigma_1 \cdot \sigma_2 \sigma_1 \sigma_3 \cdot \sigma_1 \sigma_3, \\
\nonumber  A_{2,2} &=&  \sigma_2 \sigma_1 \sigma_3 \cdot \sigma_1 \sigma_3 \cdot \sigma_1 \sigma_3 \sigma_2 \sigma_1, \\
\nonumber  A_{2,3} &=& \sigma_1 \sigma_3 \cdot \sigma_1 \sigma_3
\sigma_2 \sigma_1 \cdot \sigma_2 \sigma_1 \sigma_3.\}
\end{eqnarray}
{\rm It is easy to check (Lemma~\ref{L:tau_c_d}) that $USS(A)$ is
invariant under decycling, which for rigid braids is simply reverse
cycling.  Observe that $\tau(A_{1,j}) = A_{2,j}$ for $j=1,2,3$ and
that $\tau^2 = 1$, so that our example also illustrates the fact
that $USS(A)$ is invariant under $\tau$.

The graph $\Gamma_A$ is illustrated in Figure~\ref{F:A}.  Since all
black arrows (partial cyclings) in $\Gamma_A$ correspond to
cyclings, the only question in constructing the graph is how to
relate vertices in $A_1$ to vertices in $A_2$ by partial twisted
decyclings.  For that, we need to know the prefixes of the inverses.
By Lemma~\ref{L:initial & final factors of inverses}, for every
$X\in G$ one has $\iota(X^{-1})=\partial(\varphi(X))$, and in this
way we find that $\iota(A_{1,1}^{-1}) =
\sigma_2\sigma_1\sigma_3\sigma_2$, and in fact
$(A_{1,1})^{\sigma_2\sigma_1\sigma_3\sigma_2} = A_{2,3}$.   All
other arrows on the graph $\Gamma_A$, which can be seen in the left
sketch in Figure~\ref{F:A}, can be computed in a similar way.   In
this very simple example it happens that both $\iota(A_{i,j})$ and
$\iota(A_{i,j}^{-1})$ are {\it minimal} simple elements, for every
$i,j$. Hence, all partial cyclings are actually cyclings, and all
partial twisted decyclings are actually twisted decyclings, so they
connect $A_{i,j}$ to $\tau(\mathbf d(A_{i,j}))$. }

\end{example}
Let $\overline{A} = A^{-1}$.  The right sketch in Figure~\ref{F:A}
depicts the graph $\Gamma_{A^{-1}} = \Gamma_{\overline A}$.  There
are 6 elements in the cycling orbit of $\overline{A}$.  We have set
$\overline{A}_{1,1} = \overline{A}$ and $\overline A_{1,j} = \mathbf
c^{j-1}(\overline A), \ \ j = 2,\dots,6$.    In this case there are
3 elements in each of 2 cycling orbits in $\Gamma_A$, but 6 elements
in 1 cycling orbit for $\Gamma_{A^{-1}}$.   We will have more to say
about this and about Example $A$ (and other examples) later.   See,
in particular, Example~\ref{E:A again} below.

\begin{figure}[htpb]
\centerline{\includegraphics[scale=1.0]{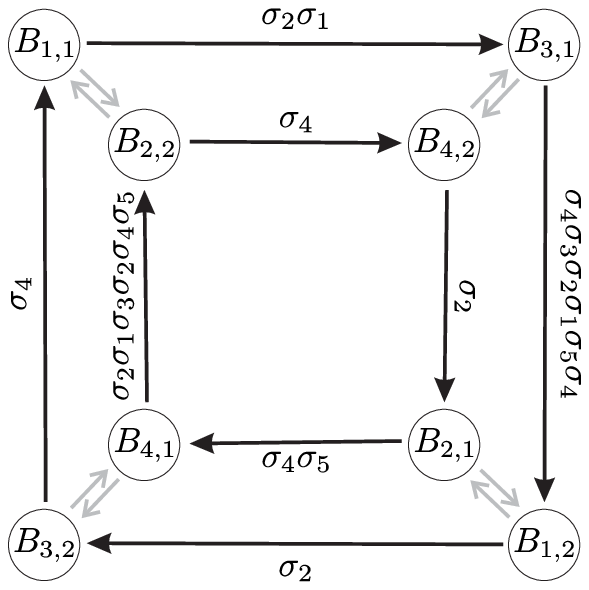}} \caption{$\Gamma_B.$
} \label{F:B}
\end{figure}

In the graph  $\Gamma_A$ all partial cyclings were actual cyclings,
but that is far from the case in general:

\begin{example}\label{E:B}
{\rm Consider the 6-braid $B =  \sigma_2 \sigma_1 \sigma_4 \sigma_3
\sigma_2 \sigma_1 \sigma_5 \sigma_4 \cdot \sigma_2 \sigma_4$  whose
graph is depicted in Figure~\ref{F:B}. Braid $B$, like braid $A$, is
PA and rigid.  A computation shows that $USS(B)$ has 4 cycling
orbits, $B_i, \ \ i = 1,2,3,4$ with 2 elements each, $B_{i,1}$ and
$B_{i,2} = \mathbf c(B_{i,1}) = \mathbf d(B_{i,1})$. Also
$\tau(B_{i,j}) = B_{i+1,j}$ for $i=1,3$ and $j=1,2$:}
\begin{eqnarray}
\nonumber  B_{1,1} &=& \sigma_2 \sigma_1 \sigma_4 \sigma_3 \sigma_2 \sigma_1 \sigma_5 \sigma_4 \cdot \sigma_2 \sigma_4,  \ \ \ B_{1,2} = \sigma_2 \sigma_4\cdot  \sigma_2 \sigma_1 \sigma_4 \sigma_3 \sigma_2 \sigma_1 \sigma_5 \sigma_4, \\
\nonumber  B_{2,1} &=&  \sigma_2 \sigma_1 \sigma_3 \sigma_4 \sigma_3 \sigma_2 \sigma_5 \sigma_4 \cdot \sigma_2 \sigma_4, \ \ \  B_{2,2} = \sigma_2 \sigma_4 \cdot  \sigma_2 \sigma_1 \sigma_3 \sigma_4 \sigma_3 \sigma_2 \sigma_5 \sigma_4 \\
 \nonumber  B_{3,1} &=&   \sigma_1 \sigma_4 \sigma_3 \sigma_2 \sigma_1 \sigma_5 \sigma_4 \cdot \sigma_2 \sigma_1 \sigma_4, \ \ \  B_{3,2} =  \sigma_2 \sigma_1 \sigma_4\cdot  \sigma_1 \sigma_4 \sigma_3 \sigma_2 \sigma_1 \sigma_5 \sigma_4 \\
 \nonumber  B_{4,1} &=&   \sigma_2 \sigma_1 \sigma_3 \sigma_2 \sigma_4 \sigma_5 \sigma_4 \cdot\sigma_2 \sigma_4 \sigma_5, \ \ \  B_{4,2} =  \sigma_2 \sigma_4 \sigma_5  \cdot  \sigma_2 \sigma_1 \sigma_3 \sigma_2 \sigma_4 \sigma_5 \sigma_4.
\end{eqnarray}
 \end{example}
In this example every black arrow is a partial cycling (see
Definition~\ref{D:partial cycling}) and the concatenation of two
consecutive black arrows corresponds to a cycling.  Hence the
product of four consecutive black arrows starting at $B_{i,j}$ is a
circuit at $B_{i,j}$, and the product of its labels is equal to
$B_{i,j}$.  The grey arrow starting at $B_{i,j}$ is
$\partial(\varphi(B_{i,j}))$ for every $i,j$. Hence, it corresponds
to a twisted decycling. We have not written down the labels for the
grey arrows because of space considerations, but will compute some
of them later, in Example~\ref{E;transport}.

\begin{figure}[htpb]
\centerline{\includegraphics[scale=1.0]{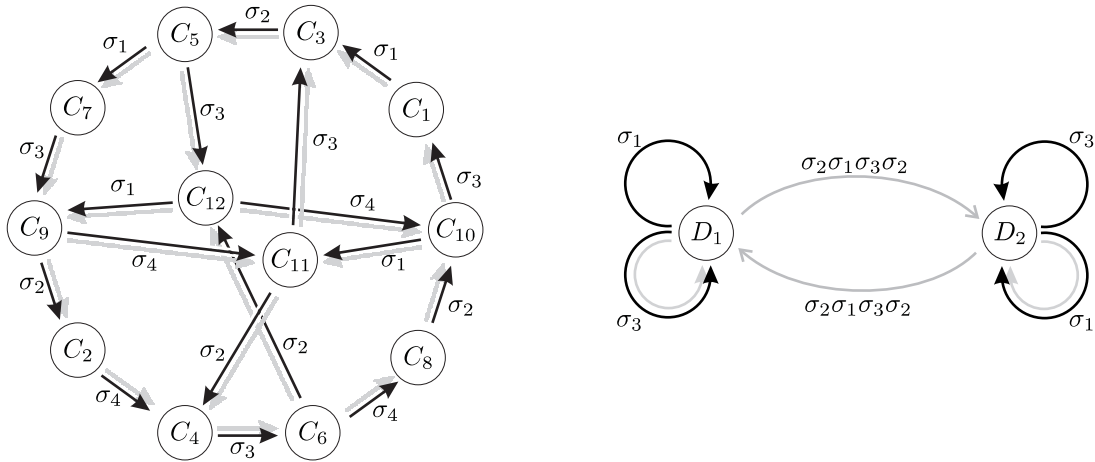}}
\caption{$\Gamma_C$ and $\Gamma_D.$ }
\label{F:C-D}
\end{figure}

In the graphs $\Gamma_A, \Gamma_{A^{-1}},$ and $\Gamma_B$ all arrows were either black or grey.  Lest the reader think that this is always the case, the left and right sketches in Figure~\ref{F:C-D} illustrate two different types of examples.
We begin with the left sketch in Figure~\ref{F:C-D}:
\begin{example}
\label{E:C}
{\rm  We consider the 5-braid $C = \sigma_ 4\sigma_ 1\sigma_ 2\sigma_ 3\sigma_ 4$.  Since $C$ is a simple braid, we have $\inf(C) = 0,  \sup(C) = 1, \ell(C) = 1$.  Its Thurston-Nielsen classification is periodic.  A computer calculation shows that $USS(C)$ has 12 elements, each of canonical length $1$. The graph $\Gamma_C$ is illustrated in the left sketch in Figure~\ref{F:C-D}.}  \end{example}
The 12 elements in $USS(C)$ are labeled $C_1,C_2,\dots, C_{12}$.  We observe that, whereas  in Example~\ref{E:A} every arrow was either black or grey, one sees that in Example~\ref{E:C} every arrow is {\it both} black {\it and} grey, that is, it is bi-colored!
Moreover, in the braid $A$ every arrow that begins at a vertex $A_{i,j}$ is labeled by either the simple element $\iota(A_{i,j})$ or the simple element $\iota(A_{i,j}^{-1})$, and so is never a proper prefix, but  in the braid $C$ every vertex is a simple element of letter length 5, and every arrow is labeled by a simple element of letter length 1 (an atom), from which it follows  that every arrow at every vertex $C_j$ is a proper prefix of  both  $C_j = \iota(C_j)$ and of $\partial(C_j)=\iota(C_j^{-1})$.  (Note that $C_j$ is simple, but $C_j^{-1}$ is not simple.)
Lest the reader think that these two extreme cases are the only ones that occur, we give one more example, in the right sketch in Figure~\ref{F:C-D}.
\begin{example}
\label{E:D}
{\rm  Consider the 4-braid $D = D_1 = \sigma_1\sigma_3\cdot \sigma_1$.  Its graph $\Gamma_D$, illustrated in the right sketch in Figure~\ref{F:C-D},  has arrows of  all three types: black, grey and
bi-colored.  The Thurston-Nielsen classification of this braid shows it to be reducible. While there is no shortage of reducible braids that are rigid, this one is not rigid.  Its graph has 2 cycling orbits, each having one element, the first that of $D_1$ and the second that of  $D_2 = \sigma_1\sigma_3\cdot \sigma_3$.  }
\end{example}

\subsection{Black and grey components of USS(X)}
\label{SS:black_grey}

The examples that we just discussed should make it clear that we need to study and try to uncover underlying structure, if possible.
Keeping in mind all of our examples,  we begin our work by asking what happens to the graph $\Gamma_X$ if one considers only black arrows, or only grey arrows.  We denote by
$\Gamma_X^{black}$ [$\Gamma_X^{grey}$] the subgraph of
$\Gamma_X$ having the same set of vertices, but just the black
[grey] arrows. The main difference between these graphs and
$\Gamma_X$ is that they are not necessarily connected. Therefore,
in general one cannot generate the whole $USS(X)$ by using only partial
cyclings [partial twisted decyclings]; both
types of conjugations are, in general, required.

We denote by ${\cal B}_1,\ldots, {\cal B}_s$ the connected
components of $\Gamma_X^{black}$ and we refer to them as the {\it black
components} of $\Gamma_X$. In the same way, we denote by ${\cal
G}_1,\ldots,{\cal G}_t$ the connected components of
$\Gamma_X^{grey}$ and we call them {\it grey components} of
$\Gamma_X$. Given $Y\in USS(X)$, we will also denote by ${\cal B}_Y$
[${\cal G}_Y$] the black [grey] component of $\Gamma_X$
containing $Y$ as a vertex.

In this subsection we will see that a black or grey component can be
computed the same way as one computes $\Gamma_X$. Starting with
any vertex $Y$, one can obtain any other vertex in ${\cal B}_Y$
[${\cal G}_Y$] by applying repeated partial cyclings [partial
twisted decyclings] to $Y$.  In other words, given any two
vertices $Y,Z \in {\cal B}_Y$ [${\cal G}_Y$] there is a
path of black [grey] arrows that goes from $Y$ to $Z$
following the sense of the arrows. In order to prove this property,
we will study the paths in $\Gamma_X$.

\begin{definition}
\label{D:path}
{\rm   An arrow in $\Gamma_X$ is characterized by its starting point
$Y\in USS(X)$ and its label $s$ which is a minimal simple element for $Y$;
the endpoint then is $Y^s$.  In order to simplify notation we denote the
arrow by its label $s$.  The {\it starting point}   and {\it endpoint} of $s_i^{-1}$ are the endpoint and starting point respectively of $s_i$.
A {\it path} in $\Gamma _X$ is a sequence
$(s_1^{e_1},\ldots,s_k^{e_k})$, possibly empty, where $s_i$ is an arrow in
$\Gamma_X$ and $e_i=\pm 1$, such that the
endpoint of $s_i^{e_i}$ is equal to the starting point of
$s_{i+1}^{e_{i+1}}$ for every $i=1,\ldots,k-1$. A path is
{\it  black} [{\it grey}], if all of its arrows are black [grey].
We say that a path $(s_1^{e_1},\ldots,s_k^{e_k})$ is {\it
oriented} if $e_i=1$ for $i=1,\ldots,k$.}
\end{definition}

Notice that every path
$(s_1^{e_1},\ldots,s_k^{e_k})$ has an associated element
$\alpha = s_1^{e_1}\cdots s_k^{e_k}$ and that different paths may have the
same associated element. Since the labels of arrows are simple elements, it
follows that if the path is oriented then $\alpha\in P$. Finally, notice
that if $X$ and $Y$ are the initial and final vertices of the above
path, that is, if $X$ is the starting point of $s_1^{e_1}$ and $Y$ is
the endpoint of $s_k^{e_k}$, then $X^\alpha=Y$.

{\bf Notation:} Given an element $Y\in G$ whose left normal form is $\Delta^p
y_1\cdots y_r$, we define $Y^{\circ}= Y\Delta^{-p}=
\tau^{-p}(y_1\cdots y_r)$.  One has the following result.

\begin{proposition}\label{P:interior->colored_path}
Let $(s_1,\ldots,s_k)$ be an oriented path in $\Gamma _X$ starting
at a vertex $Y$, and let $\alpha= s_1\cdots s_k\in P$ be its
corresponding element. Then one has:
\begin{enumerate}
\item If $\alpha \preccurlyeq Y^{\circ}$ then $(s_1,\ldots,s_k)$
is an oriented black path.

\item If $\alpha \preccurlyeq (Y^{-1})^{\circ}$ then $(s_1,\ldots,s_k)$
is an oriented grey path.

\end{enumerate}
\end{proposition}

\begin{proof}
The claim is trivial if $\alpha = 1$, so we may assume that $\alpha\not= 1$.
We will show the result by induction on $k$. Suppose that
$\alpha\preccurlyeq Y^\circ$. Notice that $\inf(Y^\circ)=0$ and $\iota(Y^
\circ)= \tau^{-p}(y_1)=\iota(Y)$. Hence, since $\alpha \preccurlyeq
Y^\circ$ it follows that $\inf(\alpha)=0$ and then
$s_1\preccurlyeq\iota(\alpha) \preccurlyeq \iota(Y^\circ)=\iota(Y)$, that is, $s_1$ is
a black arrow. In particular, the result is true if $k=1$. Suppose
that $k>1$ and that the result is true for oriented paths of length $k-1$.

We already know that $s_1$ is a black arrow. Now denote $T=Y^{s_1}$
and notice that one has $T^\circ = s_1^{-1} Y^\circ \tau^{-p}(s_1)$
since $T\in USS(X)\subset SSS(X)$. Hence, since $\alpha=s_1\cdots
s_k\preccurlyeq Y^\circ$ it follows that $s_2\cdots s_k \preccurlyeq s_1^{-1}
Y^\circ \:\preccurlyeq \:T^\circ$. By induction hypothesis $(s_2,\cdots
,s_k)$ is an oriented black path, hence $(s_1,\ldots,s_k)$ is also a
black path and the result follows.

The proof for $s\preccurlyeq (Y^{-1})^\circ$ is similar; for the induction step
note that $T^{-1}\in SSS(X^{-1})$ for any $T\in USS(X)$, whence
$(T^{-1})^\circ = s_1^{-1} (Y^{-1})^\circ \tau^{-p}(s_1)$ if $T=Y^{s_1}$.
\end{proof}

The following is a particular case of the above result, which makes
the connection with partial cyclings and partial twisted decyclings.

\begin{corollary}\label{C:initial->colored_path}
Let $(s_1,\ldots,s_k)$ be an oriented path in $\Gamma _X$ starting
at a vertex $Y$. If the associated element $s=s_1\cdots s_k$ is
simple, then one has:
\begin{enumerate}

 \item If $s\preccurlyeq\iota(Y)$ then $(s_1,\ldots,s_k)$ is an oriented black path.

 \item If $s\preccurlyeq\iota(Y^{-1})$ then $(s_1,\ldots,s_k)$ is an oriented grey path.

\end{enumerate}
\end{corollary}

\begin{proof}
 If $s$ is simple, then $s\preccurlyeq Y^\circ$ if and only if $s\preccurlyeq
 \iota(Y^\circ) = \iota(Y)$, and $s\preccurlyeq (Y^{-1})^\circ$ if and only if $s\preccurlyeq
 \iota((Y^{-1})^\circ) = \iota(Y^{-1})$. Therefore, this
 result is a direct consequence of
 Proposition~\ref{P:interior->colored_path}.
\end{proof}

\begin{example}
\label{E:T}
{\rm The converse of the above result is not true. For instance, consider
the braid $U=\sigma_1\sigma_2 \sigma_1 \sigma_3 \sigma_2\in B_4$.
Then $(\sigma_1, \sigma_3)$ is a black path in $\Gamma_U$ starting
at $U$, also $s=\sigma_1\sigma_3$ is simple, but $s$ is not a prefix
of $\iota(U)=U$.  Moreover $(\sigma_1, \sigma_3)$ is also a grey path
but $(\sigma_1\sigma_2 \sigma_1 \sigma_3 \sigma_2)
(\sigma_1\sigma_3)$ is not simple, thus $s\not\preccurlyeq \iota(U^{-1})$.
Hence this is a counterexample for the converse of the two
properties above.}
\end{example}

\begin{proposition}
Given $Y\in USS(X)$ and a black [grey] arrow $s$ in $\Gamma_X$
starting at $Y$, there exists an oriented black [grey] path
$(s_1,\ldots,s_k)$ in $\Gamma_X$ starting and ending at $Y$, such
that $s_1=s$.
\end{proposition}

\begin{proof}
Suppose that $s$ is a black arrow and let $p=\inf(Y)$. Since $Y^{Y^\circ}=
\tau^{-p}(Y)\in USS(X)$, we know by
Proposition~\ref{P:interior->colored_path} that every decomposition
of $Y^\circ$ as a product of minimal simple elements corresponds to
an oriented black path. Moreover, $s$ is a black arrow, that is, $s\preccurlyeq
\iota(Y)\preccurlyeq Y^\circ$, whence there is a decomposition of $Y^\circ$
as a product of minimal simple elements whose first factor is $s$.
Therefore, there exists a black path $(s_1,\ldots,s_t)$ in
$\Gamma_X$, going from $Y$ to $\tau^{-p}(Y)$, such that $s_1=s$.

We can apply the same reasoning to show that there is an oriented
black path from $\tau^{-(m-1)p}(Y)$ to $\tau^{-mp}(Y)$ for every
$m\geq 1$. Concatenating these paths we obtain for every $m\geq 1$ an
oriented black path from $Y$ to $\tau^{-mp}(Y)$ whose first arrow
is $s$. Since $\Delta^e$ is central for some integer $e$, it follows
that for $m=e$ there is a black path $(s_1,\ldots,s_k)$ going from
$Y$ to $\tau^{-ep}(Y)=Y$, such that $s_1=s$.

The analogous proof works for grey arrows. We just need to notice
that $Y^{(Y^{-1})^\circ}= ((Y^{-1})^{(Y^{-1})^\circ})^{-1}=
(\tau^{p+r}(Y^{-1}))^{-1} = \tau^{p+r}(Y) \in USS(X)$, where $r=\ell(Y)$.
\end{proof}

\begin{corollary}
\label{C:black and grey connected by oriented}
Given two elements $Y$ and $Z$ in a black component ${\cal B}_i$
[grey component ${\cal G}_i$] of $\Gamma_X$, there exists an
oriented black [grey] path going from $Y$ to $Z$.
\end{corollary}

\begin{proof}
Suppose that $Y$ and $Z$ belong to the same black component. Hence
one can go from $Y$ to $Z$ along a black path
$(s_1^{e_1},\ldots,s_t^{e_t})$. Suppose that $e_j=-1$ for some $j$,
where $s_j$ is a black arrow going from $U$ to $V$ (so $s_j^{-1}$
goes from $V$ to $U$). By the above result there exists an oriented
black path $(s_j,b_2,\ldots,b_k)$ going from $U$ to itself.
Therefore $(b_2,\cdots,b_k)$ is an oriented black path going from
$V$ to $U$, so one can replace $s_j^{-1}$ by $(b_2,\ldots,b_k)$ in
the above path. Applying the same procedure for every $j$ such that
$e_j=-1$, we obtain an oriented black path going from $Y$ to $Z$.
The analogous proof works for grey arrows.
\end{proof}

We establish one more property of the black and grey components:

\begin{proposition}
The set of vertices in a black component of $\Gamma_X$ is a union
of orbits under cycling.  This is not necessarily true for grey components.
\end{proposition}

\begin{proof}
To prove the first statement, we just need to show that $\mathbf c(Y)$ is a vertex of ${\cal
B}_Y$, for every $Y\in USS(X)$. But this is clear, since $\iota(Y)$
can be decomposed as a product of minimal simple elements,
$\iota(Y)=s_1\cdots s_k$. By
Corollary~\ref{C:initial->colored_path}, the path $(s_1,\cdots,
s_k)$ is black, and it goes from $Y$ to $Y^{\iota(Y)}=\mathbf c(Y)$.
Hence, $Y$ and $\mathbf c(Y)$ belong to the same black component, as
we wanted to show.

To prove that the analogous statement for grey components is false, it suffices to produce an example.
See the graph $\Gamma_B$ in Figure~\ref{F:B}. Observe that $\{B_{1,1}, B_{2,2}\}$ is not a union of orbits.
\end{proof}

Thanks to the above results, we can give an algorithm to compute the
black or the grey component of an element in its USS. It is analogous
to the corresponding algorithms computing $SSS(X)$~\cite{el-m,F-GM}
or $USS(X)$~\cite{gebhardt}. We just need to recall from
Theorem~\ref{T:convexity_USS} that given $Y\in USS(X)$ and an atom $a\in P$,
there exists a unique element $c_Y(a)$ minimal with respect to the condition
that $Y^{c_Y(a)}\in USS(Y)$ and $a\preccurlyeq c_Y(a)$.

The algorithm to compute ${\cal B}_X$ goes as follows. Starting with
$X$, conjugate it by all its minimal simple elements which are
prefixes of $\iota(X)$. For each new element $Y$ which appears in
this way, conjugate it by all its minimal simple elements which are
prefixes of $\iota(Y)$.
Keep going until no new element appears. At this point,
by Corollary~\ref{C:black and grey connected by oriented}, the black component
${\cal B}_X$ has been computed. The computation of ${\cal G}_X$ is
analogous. More precisely, the algorithms are as follows.
\bigskip\bigskip

\noindent {\sc Algorithm 1.} \hspace{1cm} {\bf Input:} $X\in
USS(X)$. \hspace{1cm} {\bf Output:} ${\cal B}_X$.

\setlength{\parskip}{0 pt}
\begin{enumerate}

 \item Set ${\cal V}=\{X\}$ and ${\cal V}'=\emptyset$.

 \item While ${\cal V}\neq {\cal V}'$ do

 \begin{enumerate}

   \item Take $Y\in {\cal V}\backslash {\cal V}'$.

   \item For every atom $a\preccurlyeq \iota(Y)$ do

   \begin{enumerate}

      \item Compute $c_Y(a)$.

      \item If $c_Y(a)$ is a minimal simple
element, set ${\cal V}={\cal V}\cup \{Y^{c_Y(a)}\}$, and store $c_Y(a)$
      as a black arrow going from $Y$ to $Y^{c_Y(a)}$.

   \end{enumerate}

   \item Set ${\cal V}'={\cal V}'\cup \{Y\}$.

 \end{enumerate}

 \item Return ${\cal V}$, together with the information on all black
 arrows.

\end{enumerate}

\bigskip\bigskip

\noindent {\sc Algorithm 2.} \hspace{1cm} {\bf Input:} $X\in
USS(X)$. \hspace{1cm} {\bf Output:} ${\cal G}_X$.

\begin{enumerate}

 \item Set ${\cal V}=\{X\}$ and ${\cal V}'=\emptyset$.

 \item While ${\cal V}\neq {\cal V}'$ do

 \begin{enumerate}

   \item Take $Y\in {\cal V}\backslash {\cal V}'$.

   \item For every atom $a$ such that $\varphi(Y)a$ is simple do

   \begin{enumerate}

      \item Compute $c_Y(a)$.

      \item If $c_Y(a)$ is a minimal simple
element, set ${\cal V}={\cal V}\cup \{Y^{c_Y(a)}\}$, and store $c_Y(a)$
      as a grey arrow going from $Y$ to $Y^{c_Y(a)}$.

   \end{enumerate}

   \item Set ${\cal V}'={\cal V}'\cup \{Y\}$.
 \end{enumerate}

 \item Return ${\cal V}$, together with the information on all grey
 arrows.

\end{enumerate}

\setlength{\parskip}{1ex plus 0.5ex minus 0.2ex}

\subsection{Intersection of black and grey components}
\label{SS:Intersection_of_components}

Now that we know how to compute the black and grey components of any
element in $USS(X)$, we will see how this can be used to solve the
conjugacy problem in $G$. Recall that the algorithm
in~\cite{gebhardt} is based on the fact that $X\in USS(X)$
and $Y\in USS(Y)$ are conjugate if and only if
$USS(X)=USS(Y)$, or equivalently, if $Y\in USS(X)$. But the
analogous statement does not hold for black or grey components, since
$\Gamma_X^{black}$ and $\Gamma_X^{grey}$ are not necessarily
connected. That is, we could have two conjugate elements $X,Y \in
USS(X)$ such that ${\cal B}_X\cap {\cal B}_Y=\emptyset$. Hence,
computing just black components (or just grey components) will not
solve the conjugacy problem, in general.

However, we will see in this section that it is actually possible to
determine whether $X\in USS(X)$ and $Y\in USS(Y)$ are conjugate just
by computing ${\cal B}_X$ and ${\cal G}_Y$, since {\it every} black
component intersects {\it every} grey component of $\Gamma_X$. Hence,
$X$ and $Y$ are conjugate if and only if ${\cal B}_X\cap {\cal
G}_Y\neq \emptyset$. This result will be achieved by showing that
every two elements in $USS(X)$ can be joined by a grey path followed
by a black path. First we need two preliminary results.

\begin{proposition}\label{P:final->colored_path}
Let $(s_1,\ldots,s_k)$ be an oriented path in $\Gamma _X$ starting
at a vertex $Y$, with $\ell(Y)>0$. If the associated element
$s=s_1\cdots s_k$ is simple, then one has:
\begin{enumerate}

\item If $\varphi(Y) s$ is left weighted then $(s_1,\ldots,s_k)$ is an oriented black path.

\item If $\varphi(Y^{-1}) s$ is left weighted then $(s_1,\ldots, s_k)$ is an oriented grey path.

\end{enumerate}
\end{proposition}

\begin{proof}
The result follows immediately from  Proposition~\ref{P:prefix} and
Corollary~\ref{C:initial->colored_path}. \end{proof}

Let us extend the above result to non-simple elements.

\begin{proposition}\label{P:final->colored_path_nonsimple}
Let $X\in G$ and $Y\in USS(X)$, with $\ell(Y)>0$. Suppose that
$\alpha\in P$ is such that $\inf(\alpha)=0$ and $Y^{\alpha}\in
USS(X)$.

If $\varphi(Y)\iota(\alpha)$ is left weighted, then $\alpha$  can be decomposed
as $\alpha=s_1\cdots s_k$, where $(s_1,\ldots,s_k)$ is an oriented
black path.

If $\varphi(Y^{-1})\iota(\alpha)$ is left weighted, then $\alpha$  can be decomposed
as $\alpha=s_1\cdots s_k$, where $(s_1,\ldots,s_k)$ is an oriented
grey path.
\end{proposition}

\begin{proof}
If $\alpha = 1$ the result follows from Propostion~\ref{P:final->colored_path}.
Let $\Delta^p y_1\cdots y_r$ be the left normal form of $Y$ and let
$\alpha_1\cdots \alpha_t$ be the left normal form of $\alpha$. Since
$\alpha_1=\alpha\wedge \Delta$, Theorem~\ref{T:convexity_USS} tells
us that $Y^{\alpha_1}\in USS(X)$. Suppose first that $\varphi(Y)
\iota(\alpha)= y_r \alpha_1$ is left weighted. Then, by
Proposition~\ref{P:final->colored_path}, $\alpha_1$ can be
decomposed as a product of black arrows.

Now consider $Z=Y^{\alpha_1}\in USS(X)\subset SSS(X)$ whose left normal form is
$\Delta^p z_1\cdots z_r$, and $\alpha'=\alpha_2\cdots \alpha_t$. The
result will follow by induction on $t$ if we show that $z_r
\alpha_2$ is left weighted. But since $y_r \alpha_1$ is left
weighted and $Z= \alpha_1^{-1} \Delta^p y_1\cdots y_r \alpha_1$, it
follows that $z_r = \beta \alpha_1$, where $y_r\succcurlyeq \beta$.
Therefore $z_r \alpha_2 = (\beta \alpha_1)\alpha_2$ is left
weighted, since so is $\alpha_1\alpha_2$. Hence,
$\varphi(Z)\iota(\alpha')= z_r \alpha_2$ is left weighted, so
induction on $t$ can be applied and the result follows.

If $\varphi(Y^{-1})\alpha$ is left weighted, we can apply a similar
reasoning to show that $\alpha=s_1\cdots s_k$ where $(s_1,\ldots,s_k)$ is an
oriented grey path.  For the induction step note that if $Z=Y^{\alpha_1}\in
USS(X)$, then $Z^{-1} = (Y^{-1})^{\alpha_1} \in SSS(X^{-1})$.
\end{proof}

Now we are able to show that every pair of elements in $USS(X)$ can
be joined by a grey path followed by a black path, and also by a
black path followed by a grey path.
For an example, see  Figure~\ref{F:B}.

\begin{theorem}\label{T:grey_black}
Let $X\in G$. Given $Y,Z\in USS(X)$, with $\ell(Y)>0$, there exists
an oriented path $(g_1,\ldots,g_s, b_1,\ldots,b_t)$ in
$\Gamma_X$ going from $Y$ to $Z$, such that $(g_1,\ldots,g_s)$
is a (possibly empty) grey path and $(b_1,\ldots,b_t)$ is a (possibly
empty) black path.  If
$Z=Y^{\alpha}$ with $\alpha\in P$, then the paths can be chosen such that
$g_1\cdots g_s b_1\cdots b_t=\alpha$.
\end{theorem}

\begin{proof}
Since $Y$ and $Z$ are conjugate, there exists some $\alpha\in G$
such that $Y^{\alpha}=Z$. Let $\Delta^m \alpha_1\cdots \alpha_t$ be
the left normal form of $\alpha$. Since some power of $\Delta$, say
$\Delta^e$, belongs to the center of $G$, it follows that
$\Delta^{ke+m} \alpha_1\cdots \alpha_t$ also conjugates $Y$ to $Z$
for every integer $k$. Hence we can assume $\inf(\alpha)\geq 0$,
that is, $\alpha\in P$.

Denote $Y^{(1)}=Y$, and let $\Delta^p y_1\cdots y_r$ be its left
normal form. Denote also $\alpha^{(1)}=\alpha = g_1\cdots g_s b_1\cdots b_t$. If $\inf(\alpha^{(1)})>0$
then $\Delta\preccurlyeq \alpha^{(1)}$, so every simple element is a prefix of
$\alpha^{(1)}$. In particular, all grey arrows for $Y^{(1)}$ are prefixes
of $\alpha^{(1)}$. Let us then choose $g_1\preccurlyeq \partial(y_r) \preccurlyeq
\alpha^{(1)}$, a grey arrow starting at $Y^{(1)}$, and denote
$Y^{(2)}=Y^{g_1}$ and $\alpha^{(2)}=g_1^{-1}\alpha^{(1)}$. We
continue this process while $\inf(\alpha^{(i)})>0$, finding grey
arrows $g_1,\ldots,g_{i}$ such that $\alpha=g_1\cdots
g_{i} \alpha^{(i+1)}$ and $Y^{(i+1)}=Y^{g_1\cdots g_{i}}$.
Since the length of possible decompositions of $\alpha$ as a product
of simple elements is finite, this process must stop, and we will
have $\alpha= g_1\cdots g_{k-1} \alpha^{(k)}$, where
$(g_1,\ldots,g_{k-1})$ is an oriented grey path and
$\inf(\alpha^{(k)})=0$. Notice that $\alpha^{(k)}$ conjugates
$Y^{(k)}=Y^{g_1\cdots g_{k-1}}$ to $Z$.

Now let $\Delta^p y_1^{(k)}\cdots y_r^{(k)}$ be the left normal form
of $Y^{(k)}$, and suppose that $y_r^{(k)} \iota(\alpha^{(k)})$ is
not left weighted. This means that $\partial(y_r^{(k)}) \wedge
\iota(\alpha^{(k)})\neq 1$.  By Theorem~\ref{T:convexity_USS}, this
element conjugates $Y^{(k)}$ to another element in $USS(X)$, hence
there is a minimal simple element $g_k\preccurlyeq \partial(
y_r^{(k)})\wedge \iota(\alpha^{(k)})$. Since $g_k\preccurlyeq
\partial(y_r^{(k)})$, one has that $g_k$ is a grey arrow and since
$g_k\preccurlyeq \iota(\alpha^{(k)})$, it follows that
$\alpha^{(k+1)}=g_k^{-1} \alpha^{(k)}\in P$.  We can continue
this process while $\varphi(Y^{(i)}) \iota(\alpha^{(i)})$ is not
left weighted, adding new arrows to our oriented grey path dividing
$\alpha$.

As above, this process must stop. We will then obtain $\alpha=
g_1\cdots g_s \alpha^{(s+1)}$, and $Y^{(s+1)}=Y^{g_1\cdots
g_s}\in USS(X)$, where $(g_1,\ldots,g_s)$ is a grey path,
$\inf(\alpha^{(s+1)})=0$ and
$\varphi(Y^{(s+1)})\iota(\alpha^{(s+1)})$ is left weighted. Since
$\alpha^{(s+1)}$ conjugates $Y^{(s+1)}$ to $Z\in USS(X)$, it follows
by Proposition~\ref{P:final->colored_path_nonsimple} that
$\alpha^{(s+1)}$ can be decomposed as a product of black arrows, as
we wanted to show.
\end{proof}

\begin{theorem}\label{T:black_grey}
Let $X\in G$. Given $Y,Z\in USS(X)$, with $\ell(Y)>0$, there exists
an oriented path $(b_1,\ldots,b_t,g_1,\ldots,g_s)$ in
$\Gamma_X$ going from $Y$ to $Z$, such that $(b_1,\ldots,b_t)$ is a
(possibly empty) black path and $(g_1,\ldots,g_s)$ is a
(possibly empty) grey path.  If
$Z=Y^{\alpha}$ with $\alpha\in P$, then the paths can be chosen such that  $b_1\cdots
b_t g_1\cdots g_s=\alpha$.
\end{theorem}

\begin{proof}
The proof is analogous to the proof of the previous result. There
exists some $\alpha\in P$ such that $Y^{\alpha}=Z$. We then
construct a sequence $\{Y^{(i)}\}_{i\geq 1}$ and
$\{\alpha^{(i)}\}_{i\geq 1}$, where $\alpha^{(1)}=\alpha$ and
$Y^{(1)}=Y$. If $\inf(\alpha^{(i)})>0$, there exists a black arrow
$b_i$ dividing $\alpha^{(i)}$, and we define
$Y^{(i+1)}=(Y^{(i)})^{b_i}$ and $\alpha^{(i+1)}=b_i^{-1}
\alpha^{(i)}$.

If $\inf(\alpha^{(i)})=0$ and $\varphi((Y^{(i)})^{-1})
\iota(\alpha^{(i)})$ is not left weighted, then there is a prefix
$\beta\preccurlyeq \iota(\alpha^{(i)})$ such that
$\varphi((Y^{(i)})^{-1})\beta$ is simple and $(Y^{(i)})^{\beta}\in
USS(X)$, for instance, one can take $\beta=
\partial(\varphi((Y^{(i)})^{-1}))\wedge \iota(\alpha^{(i)})
\preccurlyeq \partial(\varphi((Y^{(i)})^{-1})) = \iota(Y^{(i)})$
by Lemma~\ref{L:initial & final factors of inverses}.
Hence, every minimal simple element dividing $\beta$ (and thus
$\alpha^{(i)}$) is a black arrow. Therefore, if
$\varphi((Y^{(i)})^{-1})\iota(\alpha^{(i)})$ is not left weighted,
there exists a black arrow dividing $\alpha^{(i)}$ .

We keep going, finding new black arrows dividing each
$\alpha^{(j)}$, until we obtain some $\alpha^{(t+1)}$ such that
$\varphi((Y^{(t+1)})^{-1})\iota(\alpha^{(t+1)})$ is left weighted.
Then, by Proposition~\ref{P:final->colored_path_nonsimple},
$\alpha^{(t+1)}$ can be decomposed as a product of grey arrows, and
the result follows.
\end{proof}

\begin{corollary}\label{C:components_intersect}
Let $X\in G$. Given $Y,Z\in USS(X)$, one has ${\cal B}_Y \cap {\cal
G}_Z \neq \emptyset$.
\end{corollary}

\begin{proof}
This is a straightforward consequence of the previous
result. For instance, we know that there exists an oriented path
$(b_1,\ldots,b_t, g_1,\ldots,g_s)$ going from $Y$ to $Z$, such
that $(b_1,\ldots,b_t)$ is a black path and $(g_1,\ldots,g_s)$
is a grey path. If we define $V=Y^{b_1\cdots b_t}=
Z^{g_s^{-1}\cdots g_1^{-1}}$, then $V$ belongs to the same
black component as $Y$, and to the same grey component as $Z$, hence
$V\in {\cal B}_Y \cap {\cal G}_Z$.
\end{proof}

\begin{corollary}
\label{C:new criterion}
Given $X,Y\in G$, let $X'\in USS(X)$ and $Y'\in USS(Y)$. Then $X$
and $Y$ are conjugate if and only if ${\cal B}_{X'}\cap {\cal
G}_{Y'}\neq \emptyset$.
\end{corollary}

\begin{proof}
We know that $X$ and $Y$ are conjugate if and only if
$USS(X)=USS(Y)$. In this case, $X',Y'\in USS(X)$, hence
Corollary~\ref{C:components_intersect} tells us that ${\cal
B}_{X'}\cap {\cal G}_{Y'}\neq \emptyset$. Conversely, if there
exists some $V\in {\cal B}_{X'}\cap {\cal G}_{Y'}$, then $V$ is
conjugate to $X'$ and also to $Y'$. Since $X$ is conjugate to $X'$
and $Y$ is conjugate to $Y'$, it follows that $X$ and $Y$ are
conjugate.
\end{proof}

Using the above results, we will be able to obtain a new algorithm to solve the CDP/CSP in Garside groups.
See $\S$\ref{SS:a new solution to CDP/CSP in Garside groups} below.

\section{Rigid elements}
 \label{S:rigid elements}

 Rigid elements in Garside groups were studied in \cite{BGGM-I}. The structure of $\Gamma_X$ which we described in the previous
section, will be seen in this section to be particularly simple in the case of rigid elements.  This is important for the following reasons:
\be
\item It was proved in Theorems 3.21 and 3.22 of \cite{BGGM-I} that a wide class of elements in Garside groups have rigid powers. Also, that if $X$ has a rigid power, then all elements in $USS(X)$ have rigid powers. Moreover if one element in $USS(X^k)$ is rigid, then all elements in $USS(X^k)$ are rigid.
\item Corollary 3.24 of \cite{BGGM-I} asserts, in particular,  that every PA braid has a small power whose ultra summit set
consists of rigid braids, also an explicit bound for the power was found in $\S$3.5 of \cite{BGGM-I}. This allows us to study the conjugacy
problem for PA braids, using the properties of rigid
elements.
\ee

As noted earlier, the ultra summit set of a rigid
element in a Garside group often exhibits a transparently simple structure
similar to the one seen in Example~\ref{E:A}, that is,  the cycling orbits have length $\ell(X)$ or $2\ell(X)$ and the number of orbits is 2 or 1.
However,  there are examples of braids (and thus of elements in
Garside groups) whose ultra summit set is bigger than expected.
Hence, a deeper study of ultra summit sets is needed to understand
this phenomenon.   As will be seen, the combinatorics are quite complicated and the bounds we seek are well-hidden.

\subsection{Vertices and arrows in the USS graph of a rigid element}

The effect of a partial cycling on the normal form of a rigid element can be described explicitly:
\begin{lemma}\label{L:decomposition on partial cycling}
Let $X\in G$ be a rigid element, whose left normal form is $\Delta^p x_1\cdots x_r$.  Let $s\preccurlyeq \iota(X)$, such that $X^s\in SSS(X)$. Then there is a decomposition of each factor $x_i=a_ib_i$, such that $a_i,b_i$ are simple elements for $i=1,\ldots,r$, also $s=\tau^{-p}(a_1)$, and the left normal forms of $X$ and $X^s$ are
$$
X=\Delta^p (a_1b_1)(a_2b_2)\cdots (a_rb_r) \hspace{1cm} X^s=\Delta^p (b_1a_2)(b_2a_3)\cdots (b_r \tau^{-p}(a_1)).
$$
\end{lemma}
\begin{proof}
Since $s\preccurlyeq \iota(X)=\tau^{-p}(x_1)$, we can decompose $x_1=a_1b_1$, where $s=\tau^{-p}(a_1)$. By the rigidity of $X$, the left normal form of $X^s$ is precisely $\Delta^p x_1\cdots x_r s$. Hence one has $X^s = \Delta^p b_1 x_2\cdots x_r s$, where the factors $x_2\cdots x_r s$ are in left normal form. It is then well known~\cite{Michel} that the left normal form of $Y$ is given by $\Delta^p(b_1a_2)(b_2a_3)\cdots (b_r s_1)(s_2)$, where $a_ib_i=x_i$ and $s_1s_2=s$. Since we are assuming that $X^s\in SSS(X)$, that is, $\ell(X^s)=r$, one has $s_2=1$ and $s_1=s=\tau^{-p}(a_1)$.
\end{proof}

Concerning the arrows of $\Gamma_X$, the first difference
between the case of a rigid element and the general case is the
following:

\begin{lemma}\label{L:bicolored_not_rigid}
Let $X\in USS(X)$ and $\ell(X)>0$. Then $X$ is rigid if and only if
there are no bi-colored arrows starting at $X$.
\end{lemma}

\begin{proof}
If an arrow $s$ starting at $X$ is black and
grey at the same time, then it is a left divisor of both $\iota(X)$
and $\iota(X^{-1})$. This is only possible if $\iota(X)\wedge
\iota(X^{-1})\neq 1$, that is, if $X$ is not rigid.

Conversely, suppose that $X$ is not rigid. We know that
$X^{\iota(X)}=\mathbf c(X)\in USS(X)$, but also
$X^{\iota(X^{-1})}=\tau(\mathbf d(X))\in USS(X)$. Hence, by
Theorem~\ref{T:convexity_USS}, $X^{\iota(X)\wedge \iota(X^{-1})}\in
USS(X)$, so there is a minimal simple element for $X$ which is a
divisor of $\iota(X)\wedge \iota(X^{-1})$. In other words, there is
a bi-colored arrow starting at $X$.
\end{proof}

We now recall several results about the vertices in $USS(X)$, when $X$ is rigid, all proved by the authors of this paper in the manuscript \cite{BGGM-I}.

\begin{proposition}
\label{proposition:X rigid implies every Y in USS(X) rigid} Let $X$
be rigid and $\ell(X)>1$.  Then every element in the ultra summit
set of $X$ is rigid.
\end{proposition}

\begin{corollary}\label{C:USS_of_rigid}
If $X$ is rigid and $\ell(X)>1$, then $USS(X)$ is the set of rigid
conjugates of $X$.
\end{corollary}

\begin{corollary}
\label{C:USS_inverse}
If $X$ is rigid and $\ell(X)>1$, then $USS(X^{-1})$ is the set of
inverses of the elements in $USS(X)$.
\end{corollary}

\begin{corollary}
\label{C:bicolored_not_rigid}
Let $X\in USS(X)$ and $\ell(X)>1$. Then $X$ is rigid if and only if
there are no bi-colored arrows in $\Gamma_X$.
\end{corollary}

\begin{proof}
This follows directly from Lemma~\ref{L:bicolored_not_rigid} and
Proposition~\ref{proposition:X rigid implies every Y in USS(X)
rigid}.
\end{proof}

For examples, the reader may wish to go back and compare Examples~\ref{E:A} and \ref{E:B}, both of which were rigid braids, and Examples~\ref{E:C} and \ref{E:D}, which were not rigid braids.  In this regard, we note that Example~\ref{E:C} is a periodic braid (peek ahead to $\S$\ref{S:periodic elements}), whereas Example~\ref{E:D} is reducible. The additional structure in $C$ will be elucidated when we get to $\S$\ref{S:periodic elements}.

We continue our investigations of the graph of $USS(X)$ when $X$ is rigid.  We know that $X$ is rigid if and only if $X^{-1}$ is rigid.  Now we
show that the arrows starting at $X$ in $\Gamma_X$ are the same as
the arrows starting at $X^{-1}$ in $\Gamma_{X^{-1}}$.

The next lemma will be important when we study the periodic elements in the braid groups, in \cite{BGGM-III},  We will also use it in the proof of Proposition~\ref{P:min simple}.

\begin{lemma}
\label{L:periodic_length_one}
If $X\in SSS(X)$ and $\ell(X)=1$, then $X\in USS(X)$. Moreover,
$USS(X)=SSS(X)$, and it consists of the conjugates to $X$ whose
canonical length is 1.
\end{lemma}

\begin{proof}
Notice that, if $\ell(X)=1$ then $\mathbf c(X)=\tau^{-p}(X)$, where
$p=\inf(X)$.  Since some power of $\tau$ is trivial, it follows that
$\mathbf c^k(X)=X$ for some $k$. Therefore, $X\in USS(X)$.

If we apply this reasoning to every element in $SSS(X)$, it follows
that $SSS(X)\subset USS(X)$. Since $USS(X)\subset SSS(X)$ by
definition, equality holds. By definition, $SSS(X)$ consists
of the conjugates to $X$ of minimal canonical length. Since $X\in
SSS(X)$, and $\ell(X)=1$, the elements in $SSS(X)=USS(X)$ are
precisely the conjugates of $X$ of canonical length 1.
\end{proof}

\begin{proposition}
\label{P:min simple}
Let $X$ be rigid. Then  $s$ is a minimal simple element for $X$ if
and only if it is a minimal simple element for $X^{-1}$.
\end{proposition}

\begin{proof}
By definition, $\ell(X) > 0$.  By Lemma~\ref{L:periodic_length_one}, if $\ell(X)=1$ then $USS(X)=SSS(X)$ and $USS(X^{-1})=\{
  y^{-1}\mid y \in USS(X)\}$ by \cite{el-m}, whereas in the case
$\ell(X)>1$ the set $USS(X)$ is made of rigid elements by
Proposition~\ref{proposition:X rigid implies every Y in USS(X) rigid}
  and $USS(X^{-1})=\{
  y^{-1}\mid y \in USS(X)\}$ by Corollary~\ref{C:USS_inverse}.

Hence, if $s$
is a simple element such that $s^{-1}Xs\in USS(X)$, then $s^{-1}Xs$
either has canonical length 1, or is rigid. In either case,
$s^{-1}X^{-1}s$ satisfies the same property, so $s^{-1}X^{-1}s\in
USS(X^{-1})$. Hence the set of simple elements conjugating $X$ to an
element in $USS(X)$ coincides with the set of simple elements
conjugating $X^{-1}$ to an element in $USS(X^{-1})$. Since the
minimal simple elements for $X$ and for $X^{-1}$ are the
$\preccurlyeq$-minimal elements in these sets, it follows
that the sets of minimal simple elements for $X$ and for $X^{-1}$ coincide.
\end{proof}

\begin{corollary}\label{C:inverse_black grey}
If $X$ is rigid and $s$ is a minimal simple element for $X$, then
$s$ is a black [grey] arrow in $\Gamma_X$ if and only if it is
a grey [black] arrow in $\Gamma_{X^{-1}}$.
\end{corollary}

\begin{proof}
By the above result, $s$ is an arrow in $\Gamma_X$ starting at $X$
if and only if it is also an arrow in $\Gamma_{X^{-1}}$ starting at
$X^{-1}$. The statement about the color of the arrow is then given
by definition.
\end{proof}

\begin{corollary}\label{C:isomorphic_graphs}
If $X$ is rigid and $\ell(X)>1$, then $\Gamma_X$ and
$\Gamma_{X^{-1}}$ are isomorphic graphs, with the same labels but
interchanged colors.
\end{corollary}

\begin{proof}
By Corollary~\ref{C:USS_inverse} we have
$USS(X^{-1})=\{y^{-1}\mid y \in USS(X)\}$ and by
Proposition~\ref{proposition:X rigid implies every Y in USS(X)
rigid}, every element in $USS(X)$ is rigid, so
Corollary~\ref{C:inverse_black grey} can be applied to every vertex
in $\Gamma_X$. Hence, an isomorphism between $\Gamma_X$ and
$\Gamma_{X^{-1}}$ can be defined by sending every vertex $Y$ to
$Y^{-1}$, and every arrow $s$ to itself. It is an isomorphism since
every arrow $s$ in $\Gamma_X$ starting at $Y$ and ending at $Z$
(thus $s^{-1} Y s =Z$) is sent to an arrow $s$ in $\Gamma_{X^{-1}}$
starting at $Y^{-1}$, whose final vertex is precisely $s^{-1} Y^{-1}
s = Z^{-1}$.
\end{proof}

\begin{example}
\label{E:A again}
{\rm We return to the example that we gave earlier, in Figure~\ref{F:A}.  We  included all edge labels, for easy comparison between the two graphs. Recall that $\overline{A} = A^{-1}$. The reader will see that the two black cycling orbits of length 3 in $\Gamma_A$ correspond to two grey orbits in $\Gamma_{\overline A}$, whereas the black cycling orbit of length 6 in $\Gamma_{\overline A}$ corresponds to a grey orbit of length 6 in $\Gamma_A$.  The correspondence between vertex labels in $\Gamma_A$ and $\Gamma_{\overline A}$  is:
$A_{1,1}^{-1} = \overline A_{1,1}, \ \ \
A_{1,2}^{-1} = \overline A_{1,3},  \ \ \
A_{1,3}^{-1} = \overline A_{1,5},  \ \ \
A_{2,1}^{-1} = \overline A_{1,4},  \ \ \
A_{2,2}^{-1} = \overline A_{1,6}, \ \ \
A_{2,3}^{-1} = \overline A_{1,2}. $    }
\end{example}

In the case of rigid elements, there is a characterization of black
and grey arrows which does not hold in the general case.

\begin{proposition}\label{P:rigid_black_arrow}
If $X$ is rigid with $\ell(X)>0$, and $s$ is a minimal simple
element for $X$, then the following conditions are equivalent.
\begin{enumerate}

\item $s$ is a black arrow in $\Gamma_X$ (that is, $s\preccurlyeq \iota(X)$).

\item $\varphi(X^{-1})s$ is simple.

\item $\varphi(X)s$ is left weighted.

\item $s\wedge \iota(X^{-1})=1$.

\end{enumerate}
The equivalence of conditions 1 and 2 and the equivalence of
conditions 3 and 4 hold for arbitrary $X\in USS(X)$.
\end{proposition}

\begin{proof}
For arbitrary $X\in USS(X)$, the product $\varphi(X)s$ is left weighted
if and only if $s\wedge \partial(\varphi(X))=1$. By
Lemma~\ref{L:initial & final factors of inverses}, we know that
$\partial(\varphi(X)) = \iota(X^{-1})$, hence conditions 3 and 4 are
equivalent.  Moreover, $\iota(X)=\partial(\varphi(X^{-1}))$, so $s$
is black if and only if $s\preccurlyeq \partial(\varphi(X^{-1}))$ or, in
other words, if and only if $\varphi(X^{-1})s$ is simple. Hence,
conditions 1 and 2 are equivalent.

Suppose that $X$ is rigid with $\ell(X)>0$, we know by
Theorem~\ref{T:minimal_simple_elements} that $\varphi(X)s$ is either
simple or left weighted. By Lemma~\ref{L:bicolored_not_rigid} we
know that $s$ is black if and only if it is not grey, which happens,
since $\iota(X^{-1})=\partial(\varphi(X))$, if and only if
$\varphi(X)s$ is not simple. Therefore, conditions 1 and 3 are
equivalent.
\end{proof}

\begin{proposition}\label{P:rigid_grey_arrow}
If $X$ is rigid with $\ell(X)>0$, and $s$ is a minimal simple
element for $X$, then the following conditions are equivalent.
\begin{enumerate}

\item $s$ is a grey arrow in $\Gamma_X$ (that is, $s\preccurlyeq \iota(X^{-1})$).

\item $\varphi(X)s$ is simple.

\item $\varphi(X^{-1})s$ is left-weighted.

\item $s\wedge \iota(X)=1$.

\end{enumerate}
The equivalence of conditions 1 and 2 and the equivalence of
conditions 3 and 4 hold for arbitrary $X\in USS(X)$.
\end{proposition}

\begin{proof}
The equivalence of conditions 1 and 2 and the equivalence of
conditions 3 and 4 can be shown as in the proof of
Proposition~\ref{P:rigid_black_arrow} with $X$ replaced by $X^{-1}$.

The rest is a consequence of Proposition~\ref{P:rigid_black_arrow},
together with Corollary~\ref{C:inverse_black grey}.
\end{proof}

\subsection{Paths in the USS graph of a rigid element}

In the general case, we saw how two elements in $USS(X)$ are always
joined by an oriented black path followed by an oriented grey path,
and vice versa. In the rigid case, we will see that we can say much
more. First, we can characterize black or grey oriented paths
according to their associated conjugating element.

\begin{lemma}\label{L:rigid_black_path}
Let $X\in G$ be rigid, and $\ell(X)>1$. Let $(s_1,\ldots,s_k)$ be an
oriented path in $\Gamma _X$ starting at a vertex $Y$, and let
$\alpha=s_1\cdots s_k$. The following properties are equivalent:
\begin{enumerate}

 \item The path $(s_1,\ldots, s_k)$ is black.

 \item $\alpha\preccurlyeq (Y^m)^\circ$ for some positive integer $m$.

 \item $\alpha\wedge \iota(Y^{-1})=1$.

\end{enumerate}
\end{lemma}

\begin{proof}
We will show the result by induction on $k$. The case $k=0$ is trivial.  If $k=1$ then
$\alpha=s_1$ is simple, hence the result follows from the definition
of a black arrow, since in this case $s_1\preccurlyeq (Y^m)^\circ$ if and
only if $s_1\preccurlyeq \iota((Y^m)^\circ)= \iota(Y^m)$. Note that $Y$ is rigid
by Proposition~\ref{proposition:X rigid implies every Y in USS(X)
rigid}.  Therefore  $\iota(Y^m)=\iota(Y)$, hence $s_1\preccurlyeq (Y^m)^\circ$ if and
only if $s_1\preccurlyeq \iota(Y)$, that is, if and only if $s_1$ is a
black arrow. Hence Properties 1 and 2 are equivalent if $k=1$. In the case
$\alpha=s_1$, Property~3 just means that $s_1$ is not a grey arrow, which
by Lemma~\ref{L:bicolored_not_rigid} and the rigidity of $Y$,
is equivalent to Property~1.

Suppose that $k>1$ and that the result is true for $k-1$. Let $\Delta^p
y_1\cdots y_r$ be the left normal form of $Y$ and define
$Z=s_1^{-1}Ys_1$. Notice that $(s_2,\ldots,s_k)$ is an oriented
path starting at $Z$. Notice also that, since $Z\in USS(X)$, one has
$Z^\circ= s_1^{-1} Y^\circ \tau^{-p}(s_1)$. Moreover, we know by
Proposition~\ref{proposition:X rigid implies every Y in USS(X)
rigid} that $Y$ and $Z$ are rigid, so $Y^m$ and $Z^m$ are also
rigid, hence $Y^m,Z^m\in USS(X^m)$ for every $m>0$. Therefore
$(Z^m)^\circ = (s_1^{-1}Y^m s_1)^\circ = s_1^{-1} (Y^m)^\circ
\tau^{-pm}(s_1)$ for every $m>0$.

By induction hypothesis, $(s_2,\ldots,s_k)$ is a black path if and
only if $s_2\cdots s_k \preccurlyeq (Z^m)^\circ$ for some positive $m$.
Therefore, if $(s_1\ldots ,s_k)$ is a black path, then
$\alpha=s_1s_2\cdots s_k \preccurlyeq s_1 (Z^m)^\circ = s_1 s_1^{-1}
(Y^m)^\circ \tau^{-pm}(s_1) =(Y^m)^\circ \tau^{-pm}(s_1)$. But $s_1$
is a black arrow, hence $s_1\preccurlyeq \iota(Y)=\tau^{-p}(y_1)$, so
$(Y^m)^\circ \tau^{-pm}(s_1) \preccurlyeq (Y^m)^\circ \tau^{-p(m+1)}(y_1)
\preccurlyeq (Y^{m+1})^\circ$, where the last claim is true due to the
rigidity of $Y$. Therefore, if $(s_1,\ldots, s_k)$ is a black path,
then $\alpha\preccurlyeq (Y^{m+1})^\circ$ for some positive $m$, so Property~1
implies Property~2.

Conversely, suppose that $\alpha\preccurlyeq (Y^m)^\circ$ for some positive
integer $m$. Then $s_1\preccurlyeq \alpha\preccurlyeq (Y^m)^\circ$, and since $s_1$ is
simple $s_1\preccurlyeq \iota((Y^m)^\circ)=\iota(Y^m)=\iota(Y)$, whence
$s_1$ is a black arrow. Moreover, $\alpha=s_1\cdots s_k \preccurlyeq
(Y^m)^\circ \preccurlyeq (Y^m)^\circ \tau^{-mp}(s_1)$, so $s_2\cdots s_k
\preccurlyeq s_1^{-1} (Y^m)^\circ \tau^{-mp}(s_1) = (Z^m)^\circ$. By
induction hypothesis this means that $(s_2,\ldots, s_k)$ is a black
path, hence $(s_1,\ldots,s_k)$ is a black path, that is, Property~2
implies Property~1. The first two properties are thus equivalent.

Now suppose again that $\alpha\preccurlyeq (Y^m)^\circ$ for some positive $m$.
In particular, one has $\inf(\alpha)=0$, so it follows that
$\iota(\alpha)\preccurlyeq \iota((Y^m)^\circ) = \iota(Y^m) = \iota(Y)$,
where the last equality follows from the rigidity of $Y$. Since
$Y$ is rigid, $\iota(Y)\wedge \iota(Y^{-1})=1$, hence
$\iota(\alpha)\wedge \iota(Y^{-1})=1$ and therefore $\alpha\wedge
\iota(Y^{-1})=1$. This shows that Property~2 implies Property~3.

Finally, suppose that $\alpha\wedge \iota(Y^{-1})=1$. This implies
that $\inf(\alpha)=0$, so the left normal form of $\alpha$ has the
form $\alpha_1\cdots \alpha_t$. Since $\alpha\wedge \iota(Y^{-1})=
\alpha \wedge \partial(y_r) = 1$, the left normal form of
$Y\alpha$ is precisely $\Delta^p y_1\cdots y_r \alpha_1\cdots
\alpha_t$. Moreover, since $Y$ is rigid, it follows that for every
$m>0$ the left normal form of $Y^m \alpha$ is
$$
Y^m \alpha = \Delta^{pm} \tau^{p(m-1)}(y_1)\cdots \tau^{p(m-1)}(y_r)
\tau^{p(m-2)}(y_1) \cdots \tau^{p(m-2)}(y_r) \cdots \cdots y_1\cdots
y_r \alpha_1 \cdots \alpha_t.
$$
That is, there is no contribution of $\alpha$ to the first $mr$
non-$\Delta$ factors of the left normal form of $Y^m \alpha$. However,
$\alpha^{-1}Y^m \alpha =Z^m \in USS(X^m)$, whence $\alpha \preccurlyeq
(Y^m\alpha)^\circ$ for every $m>0$. This means that, for $m$ big
enough (such that $t\leq mr$) one has $\alpha \preccurlyeq (Y^m)^\circ$.
Hence Property~3 implies Property~2, thus the last two properties
are also equivalent.
\end{proof}

Here is the analogous result for grey paths.

\begin{lemma}\label{L:rigid_grey_path}
Let $X\in G$ be rigid, and $\ell(X)>1$. Let $(s_1,\ldots,s_k)$ be an
oriented path in $\Gamma _X$ starting at a vertex $Y$, and let
$\alpha = s_1\cdots s_k$. Then the following properties are
equivalent:
\begin{enumerate}

 \item The path $(s_1,\ldots,s_k)$ is grey.

 \item $\alpha \preccurlyeq (Y^{-m})^\circ$ for some positive integer $m$.

 \item $\alpha \wedge \iota(Y) =1$.

\end{enumerate}
\end{lemma}

\begin{proof}
By Corollary~\ref{C:isomorphic_graphs} this is analogous to the
above result, considering the graph $\Gamma_{X^{-1}}$ instead of
$\Gamma _X$.
\end{proof}

In the particular case where $\alpha$ above is simple, the
characterization of black or grey paths can be made more precise.

\begin{corollary}\label{C:rigid_simple_black_path}
Let $X\in G$ be rigid and $\ell(X)>1$. Let $(s_1,\ldots,s_k)$ be an
oriented path in $\Gamma _X$ starting at a vertex $Y$ and let $s =
s_1\cdots s_k$. If $s$ is simple, the following properties are
equivalent:
\begin{enumerate}

 \item The path $(s_1,\ldots, s_k)$ is black.

 \item $s \preccurlyeq \iota(Y)$.

 \item $ s \wedge \iota(Y^{-1})=1$.

 \item $ \varphi(Y^{-1})s$ is simple.

\end{enumerate}
\end{corollary}

\begin{proof}
The equivalence of Properties 1-3 is a consequence of
Lemma~\ref{L:rigid_black_path}, together with the fact that, if $s$
is simple and $Y$ is rigid, $s\preccurlyeq (Y^m)^\circ$ if and only if
$s\preccurlyeq \iota(Y^m)^\circ=\iota(Y)$.  The equivalence of Properties 2
and 4 follows from $\iota(Y)=\partial({\varphi(Y^{-1}))}$.
\end{proof}

\begin{corollary}\label{C:rigid_simple_grey_path}
Let $X\in G$ be rigid, and $\ell(X)>1$. Let $(s_1,\ldots,s_k)$ be an
oriented path in $\Gamma _X$ starting at a vertex $Y$, and let $s =
s_1\cdots s_k$. If $s$ is simple, the following properties are
equivalent:
\begin{enumerate}

 \item The path $(s_1,\ldots, s_k)$ is grey.

 \item $s \preccurlyeq \iota(Y^{-1})$.

 \item $ s \wedge \iota(Y)=1$.

 \item $ \varphi(Y)s$ is simple.

\end{enumerate}
\end{corollary}

\begin{proof}
This result is the same as the above one, but applied to
$\Gamma_{X^{-1}}$ instead of $\Gamma_X$, thanks to
Corollary~\ref{C:isomorphic_graphs}.
\end{proof}

\begin{remark}
Using the same kind of proofs as above, one can show that
Corollaries~\ref{C:rigid_simple_black_path} and
\ref{C:rigid_simple_grey_path} are also true if $Y$ is rigid and
$\ell(Y)=1$. But this cannot be generalized to paths whose
associated element is not simple, so we will keep considering rigid
elements of canonical length greater than 1.
\end{remark}

Now recall that we showed, in the general case, that every oriented
path could be transformed into the concatenation of a black path and
a grey path, whose associated element was the same as the original
one. We will now see that, if $X$ is rigid and $\ell(X)>1$,  the
elements associated to those oriented black and grey paths are, in
fact, determined by the associated element of the whole path.


\begin{lemma}\label{L:rigid_black_grey_path}
Let $X$ be a rigid element with $\ell(X)=r>1$. Let $(b_1,\ldots,b_s,
g _1,\ldots,g _t)$ be an oriented path in $\Gamma _X$ starting
at $Y$, such that $(b_1,\ldots,b_s)$ is a black path and $(g
_1,\ldots,g _t)$ is a grey one. Let $b=b_1\cdots b_s$,\; $g
=g _1\cdots g _t$, and $\alpha =bg $.  Let $m$ be an
integer big enough such that $\ell(b)\leq mr$. Then, $b= (Y^m)^\circ
\wedge \alpha $.  In particular, $b$ and $g $
are uniquely determined by $\alpha$.
\end{lemma}

\begin{proof}
Since $(b_1,\ldots,b_s)$ is a black path,
Lemma~\ref{L:rigid_black_path} tells us that $b\preccurlyeq (Y^m)^\circ$
for some $m$ big enough. Actually, by the rigidity of $Y$, it
suffices to take $m$ such that $\ell(b)\leq \ell((Y^m)^\circ) = mr$.
Since we also have $b\preccurlyeq bg =\alpha $, it follows that $b\preccurlyeq
(Y^m)^\circ \wedge \alpha $.

Now let $Z=b^{-1}Yb$. since $Z$ is also rigid, both $Y^m$ and $Z^m$
belong to $USS(X^m)$. Hence $(Z^m)^\circ
=b^{-1}(Y^m)^\circ \tau^{-pm}(b)$. Moreover, $\iota((Z^m)^\circ) =
\iota(Z)$. However, $(g_1,\ldots,g_t)$ is a grey path starting at
$Z$, hence by Lemma~\ref{L:rigid_grey_path} $g\wedge \iota(Z)=
g \wedge \iota((Z^m)^\circ)= 1$. Since $\inf((Z^m)^\circ)=0$,
this implies that $g \wedge (Z^m)^\circ =1$, that is, $g
\wedge (b^{-1}(Y^m)^\circ \tau^{-pm}(b)) =1 $. As $b\preccurlyeq
(Y^m)^\circ$, we in particular have $g\wedge b^{-1}(Y^m)^\circ =1$.
Multiplying on the left by $b$, one gets $\alpha \wedge (Y^m)^\circ
=b$, as we wanted to show.
\end{proof}

The analogous result holds for the concatenation of a grey path
followed by a black path.

\begin{lemma}\label{L:rigid_grey_black_path}
Let $X$ be a rigid element with $\ell(X)=r>1$. Let $(g
_1,\ldots,g _t,b_1,\ldots,b_s)$ be an oriented path in $\Gamma
_X$ starting at $Y$, such that $(g _1,\ldots,g _t)$ is a grey
path and $(b_1,\ldots,b_s)$ is a black one. Let $g =g _1\cdots
g _t$, \; $b=b_1\cdots b_s$, and $\alpha =g b$.  Let $m$ be an
integer big enough such that $\ell(g)\leq mr$.  Then,
$g= (Y^{-m})^\circ \wedge \alpha$. In particular, $g$ and
$b$ are uniquely determined by $\alpha$.
\end{lemma}

\begin{proof}
Corollary~\ref{C:isomorphic_graphs} implies that this result is
equivalent to the previous one.
\end{proof}

In the case that $\alpha$ above is simple, one can simplify the
characterization of the elements $b$ and $g$. We will state the
two analogous results since we will use them later.

\begin{corollary}\label{C:rigid_simple_black_grey_path}
Let $X$ be a rigid element with $\ell(X)>1$. Let $(b_1,\ldots,b_s,
g _1,\ldots,g _t)$ be an oriented path in $\Gamma _X$ starting
at $Y$, such that $(b_1,\ldots,b_s)$ is a black path and $(g
_1,\ldots,g _t)$ is a grey one. Let $b=b_1\cdots b_s$,\; $g
=g _1\cdots g _t$ and $s =bg $.  If $s$ is simple, then
$b= \iota(Y) \wedge s$.  In particular, $b$ and $g $ are uniquely
determined by $s$.
\end{corollary}

\begin{proof}
This follows from Lemma~\ref{L:rigid_black_grey_path}, since by the
simplicity of $s$ and the rigidity of $Y$, one has $b= (Y^m)^\circ
\wedge s = \Delta \wedge (Y^m)^\circ \wedge s  = \iota((Y^m)^\circ)
\wedge s = \iota(Y^m) \wedge s = \iota(Y) \wedge s$.
\end{proof}

\begin{corollary}\label{C:rigid_simple_grey_black_path}
Let $X$ be a rigid element with $\ell(X)>1$. Let $(g
_1,\ldots,g _t,b_1,\ldots,b_s)$ be an oriented path in $\Gamma
_X$ starting at $Y$, such that $(g _1,\ldots,g _t)$ is a grey
path and $(b_1,\ldots,b_s)$ is a black one. Let $g =g _1\cdots
g _t$, \; $b=b_1\cdots b_s$ and $s =g b$. If $s$ is simple,
then $g=\partial({\varphi(Y)}) \wedge s$, that is, the left
normal form of $\varphi(Y) s$ is $(\varphi(Y)g)\:b$.  In
particular, $g $ and $b$ are uniquely determined by $s$.
\end{corollary}

\begin{proof}
This follows by applying Corollary~\ref{C:rigid_simple_black_grey_path}
to $\Gamma_{X^{-1}}$.
\end{proof}

\noindent {\bf Remark:} In the above results, the elements $b$ and
$g $ are determined by $\alpha $ (or $s$). But this
does {\bf not} mean that the paths $(b_1,\ldots,b_s)$ and $(g
_1,\ldots,g _t)$ are determined by $\alpha $, since there could
exist distinct paths $(b_1,\ldots,b_s)$ and $(b_1',\ldots,b_u')$
such that $b=b_1\cdots b_s=b_1'\cdots b_u'$.  The same could
happen for $g $.

The situation is more interesting if each of the above paths
consists only of one arrow. Consider then the following situation.
Let $X=\Delta ^p x_1\cdots x_r$, \; $Y=\Delta ^p y_1\cdots y_r$  and
$Z=\Delta ^p z_1\cdots z_r$ be three vertices of $\Gamma_X$, where
$X$ is rigid and $\ell(X)>1$. Suppose that $X$ and $Y$ are joined by
a black arrow $b$, while $Y$ and $Z$ are joined by a grey arrow $g
$. That is,
$$
\begin{CD}
  X   @>b>>  Y \\
   @.  @VVg V \\
  @.  Z
\end{CD}
$$

\begin{proposition}\label{P:b_rho_->_rho_b}
In the above situation, the product $bg$ is simple, and there
exist a unique grey arrow $g'$ and a unique black arrow $b'$ such
that $bg= g'b'$.

$$
\begin{CD}
  X   @>b>>  Y \\
  @Vg 'VV  @VVg V \\
  T @>b'>>  Z
\end{CD}
$$

Specifically, $g'=\partial(\varphi(X))\wedge bg$ and
$b'=(g')^{-1}bg$.
\end{proposition}

\begin{proof}
Since $b\preccurlyeq \iota(X)$ ($b$ is a black arrow), one can decompose
$\iota(X)=bc$ so that $Y=\Delta^p \tau^p(c) x_2\cdots x_r b$, where
$x_2\cdots x_rb$ is in left normal form as written, by the rigidity
of $X$. Hence, when we compute the left normal form of $Y$, we see
that $y_r= \omega b$ for some $\omega$ such that $x_r\succcurlyeq \omega$.
However, $g$ is a grey arrow for $Y$, hence $y_r g = \omega bg$
is simple, thus $bg$ is also simple.

Consider $b'$ and $g'$ as defined above. We know that,
since $(b,g)$ is an oriented path going from $X$ to $Z$, there
exists another path $(g _1,\ldots,g _t,b_1,\ldots,b_s)$ also
going from $X$ to $Z$, which is the concatenation of a grey path and
and black path, and such that $bg= g_1\cdots g_t b_1\cdots
b_s$. Moreover, by Corollary~\ref{C:rigid_simple_grey_black_path},
$g'=g_1\cdots g_t$ and $b'=b_1\cdots b_s$. We just need to
show that $t=s=1$, that is, $g'=g_1$ and $b'=b_1$.

Firstly, we cannot have $t=0$ or $s=0$, since otherwise we would have
an oriented black path, or an oriented grey path, whose associated
element is $bg$. But this is not possible, since in that case
$(b,g)$ would be either a black path or a grey path, by
Lemmas~\ref{L:rigid_black_path} or \ref{L:rigid_grey_path}, and we
know this is not true.

Suppose that $t>1$. Then $(g_t,b_1,\ldots,b_s)$ is a path which
has black and grey arrows, that can be transformed into a path
$(b_1',\ldots, b_u',g_1',\ldots,g_v')$ where every $b_i'$ is
black, every $g_i'$ is grey, and $v> 0$ by the above reasoning.
Hence $(g_1,\ldots, g_{t-1}, b_1',\ldots,
b_u',g_1',\ldots,g_v')$ is a path going from $X$ to $Z$. Now
$(g_1,\ldots, g_{t-1}, b_1',\ldots, b_u')$ can also be replaced
by a path $(b_1'',\ldots, b_\mu'',g_1'',\ldots,g_\nu'')$, where
every $b_i''$ is black, every $g_i''$ is grey, and $\nu> 0$.
Therefore, if $t>1$ we could obtain a path $(b_1'',\cdots
b_\mu'',g_1'',\ldots,g_\nu'',g_1',\ldots,g_v')$, whose
associated element is $bg$, and $\nu+v>1$. But by
Lemma~\ref{L:rigid_black_grey_path},
$g=g_1^{\prime\prime}\cdots g_\nu^{\prime\prime} g_1'\cdots g_v'$, that is,
the grey arrow $g$ could be decomposed into a product of more than
one grey arrow. This contradicts the fact that $g$ is a minimal
simple element, since in that case $g_1''$ would be a proper
prefix of $g$. Therefore $t=1$.

We can use the same reasoning to show that $s=1$, so $(g',b')$ is
an oriented path, where $g'$ is a grey arrow and $b'$ is a black
arrow. The uniqueness of $g'$ and $b'$ follows from
Lemma~\ref{L:rigid_grey_black_path}.
\end{proof}

We have then shown how to transform any black-grey path $(b,g)$
into a grey-black path $(g',b')$ such that $bg =g 'b'$.
Moreover, by Corollaries~\ref{C:rigid_simple_black_grey_path} and
\ref{C:rigid_simple_grey_black_path}, these two decompositions $bg
=g 'b'$ are the only two ways to decompose $bg$ into a product
of minimal simple elements.

By symmetry, if we start with a grey-black path $(g, b)$, we can
transform it into a black-grey path $(b',g ')$, where $g
b=b'g '$.

\begin{proposition}
Let $X\in G$ be rigid and $\ell(X)>1$. Let $(g ,b)$ be a path in
$\Gamma _X$ starting at $X$, where $g $ is a grey arrow and $b$ is
a black one. If we define $b'=g b\wedge \iota(X)$ and
$g'=(b')^{-1}g b$, then $(b',g ')$ is a path in $\Gamma _X$
such that $b'$ is black, $g '$ is grey, and $g b=b'g '$.  Moreover, $b'$ and $g'$ are the unique such arrows.
\end{proposition}

\begin{proof}
The result follows by applying Proposition~\ref{P:b_rho_->_rho_b} to
the graph $\Gamma _{X^{-1}}$ and using
$\partial({\varphi(X^{-1})})=\iota(X)$.
\end{proof}

We have thus shown that, if $X$ and $Y$ are two rigid elements
connected by a black arrow, then every grey arrow for $Y$ is related
to a grey arrow for $X$.   However, we do not yet know whether every
grey arrow of $X$ is related to one of $Y$.
Similarly, if $X$ and $Y$ are two rigid elements
connected by a grey arrow, the question arises of whether every
black arrow of $X$ is related to one of $Y$.

Using the language of Proposition~\ref{P:b_rho_->_rho_b}, we need to
show that given a black arrow $b$ and a grey arrow $g'$ for $X$, a
grey arrow $g$ for $X^b$ and a black arrow $b'$ for $X^{g'}$
exist such that $(X^b)^{g}=(X^{g'})^{b'}$. In the next
subsection, we will see that this is true.

As a consequence, we will show that every two elements in the same
black component have the same number of grey arrows and that  every two
elements in the same grey component have the same number of black
arrows.

\subsection{Partial transport}

Let $X=\Delta^p x_1\cdots x_r\in G$ be in left normal form, such that
$X\in SSS(X)$.
We start this subsection by recalling the
definition of the transport of a simple element, given
in~\cite{gebhardt}. Let $s$ be a simple element such that $s^{-1}Xs
= Y \in SSS(X)$. We can write this by
$X\stackrel{s}{\longrightarrow} Y$.
We know from \cite{el-m} that $\mathbf c(X)$ and $\mathbf
c(Y)$ also belong to $SSS(X)$. Notice that
$X\stackrel{\iota(X)}{\longrightarrow} \mathbf c(X)$ and
$Y\stackrel{\iota(Y)}{\longrightarrow} \mathbf c(Y)$.
In~\cite{gebhardt} it is shown that, in this situation,
$\mathbf c(X)$ and $\mathbf c(Y)$ are also conjugate by a simple element, $s^{(1)}$, where $s^{(1)}$ is
defined as the element making the following diagram commutative in the
sense explained below:
$$
\begin{CD}
 X   @>\iota(X)>> \mathbf c(X) \\
@VsVV   @VVs^{(1)}V \\
 Y @>\iota(Y)>> \mathbf c(Y)
\end{CD}
$$
More precisely, $s^{(1)}= \iota(X)^{-1} \:s \:\iota(Y)$.  The
nontrivial fact shown in~\cite{gebhardt} is that $s^{(1)}$ is
simple.  The simple element $s^{(1)}$ is the {\it transport} of $s$.

\noindent {\bf Remark:} When we deal with a diagram such as the previous
one, in which the arrows represent conjugating elements, then by saying
that the diagram is {\it commutative} we mean the following: for every
two paths in the diagram whose starting and ending vertices
coincide, the product of the arrows forming every path (with the
corresponding signs) are the same. For instance, in the above
diagram this is equivalent to $\iota (X)\: s^{(1)} = s
\:\iota (Y)$, that is, $s^{(1)}=\iota(X)^{-1} \:s
\:\iota(Y)$.

We can then define the transport of $s^{(1)}$, denoted $s^{(2)}$,
and so on. In general, the transport of $s^{(i-1)}$ is denoted
$s^{(i)}$ and the above definition tells us that
$$
    s^{(i)} =  \iota(\mathbf c^{i-1}(X))^{-1} \cdots \iota(\mathbf c(X))^{-1} \iota(X)^{-1}\; s \;
              \iota(Y) \iota(\mathbf c(Y)) \cdots \iota(\mathbf c^{i-1}(Y)).
$$
We can easily see that, if $X\in USS(X)$ is rigid and $\ell(X)>1$, the
situation is much simpler. In this case, if the left normal form of
$X$ is $\Delta^p x_1\cdots x_r$, that is $\iota(X)=\tau^{-p}(x_1)$,
rigidity of $\mathbf c^k(X)$ for all $k\ge 0$ implies
$\iota(\mathbf c(X))=\tau^{-p}(x_2)$, $\iota(\mathbf c^2(X))=\tau^{-p}(x_3)$ and so
on. We also know that $Y= \Delta^p y_1\cdots y_r$ must be rigid, so
the same formulae hold for $Y$. Hence, for $i=1,\ldots,r$ we have
$$
  s^{(i)}= \tau^{-p}(x_i^{-1}\cdots x_1^{-1}) \; s \; \tau^{-p}(y_1\cdots y_i).
$$
Eventually, one has
$$
   s^{(r)}= \tau^{-p}(x_r^{-1} \cdots x_1^{-1}) \; s \; \tau^{-p}(y_1\cdots y_r) =
    \Delta^{p} X^{-1} s \; Y \Delta^{-p}
$$
$$
    = \Delta^{p}X^{-1} s (s^{-1} X s) \Delta^{-p}=  \Delta^p s \: \Delta^{-p} =  \tau^{-p}(s).
$$

Hence, if $m$ is an integer such that $pm$ is a multiple of $e$
(where $\Delta^e$ is central), one has $ s^{(rm)}= \tau^{pm}(s) =
s$. For instance, in braid groups with the Garside structure given
by the Artin generators, one has $e=2$, so we always have
$s^{(2r)}=s$.  If $p$ is even, $s^{(r)}=s$ holds.

There is another way of looking at the transports of $s$ in the rigid
case. By definition, we have $s^{(1)}=\tau^{-p}(x_1)^{-1}s\:
\tau^{-p}(y_1)$, whence $y_1=\tau^p(s^{-1})\: x_1 \:\tau^p(s^{(1)})$.
But since $X$ and $Y$ are rigid, we have $\iota(\mathbf
c(X))=\tau^{-p}(x_2)$ and $\iota(\mathbf c(Y))=\tau^{-p}(y_2)$,
whence, $y_2= \tau^p(s^{(1)})^{-1}\: x_2 \:\tau^p(s^{(2)})$.  In
general, $y_i= \tau^p(s^{(i-1)})^{-1} x_i \tau^p(s^{(i)})$ for
$i=1,\ldots,r$.  Recalling $s^{(r)}=\tau^{-p}(s)$, the
normal form of $Y$ is
$$
   Y= \Delta ^p y_1 y_2\cdots y_r  =
   \Delta ^p \left( \tau^p(s)^{-1} x_1 \: \tau^p(s^{(1)})\right) \left( \tau^p(s^{(1)})^{-1} x_2 \: \tau^p(s^{(2)})\right) \cdots
   \left( \tau^p(s^{(r-1)})^{-1} x_r \:s\right).
$$
Hence, {\bf in the rigid case}, the first $r$ transports of $s$ are
precisely the simple elements that tell us how to relate the left
normal forms of $X$ and $Y$ (conjugated by $\Delta^{-p}$).

By Corollary 2.7 of [14] the transport is a bijection  from the minimal
simple elements for $X$ to the minimal simple elements for $\mathbf c(X)$.
Hence the map $\phi$ from the graph $\Gamma_X$ to itself defined by sending
every vertex $Y$ to $\mathbf c(Y)$ and every arrow $s$ to $s^{(1)}$ is an
automorphism of $\Gamma_X$.
It is also shown in~\cite{gebhardt} that $s\preccurlyeq t$ implies
$s^{(1)}\preccurlyeq t^{(1)}$ and that $(\iota(Y))^{(1)} =
\iota(\mathbf{c}(Y))$.  By definition this means that $\phi$ sends black
arrows to black arrows.  It is elementary to check that the image of a grey
arrow is a grey arrow.  Hence $\phi$ is an automorphism of $\Gamma_X$
preserving the colors of arrows.

\begin{example}\label{E;transport} {\rm Recall that the braids we called $A$ and $B$, with $\Gamma_A$ illustrated in the left sketch in Figure~\ref{F:A} and $\Gamma_B$ illustrated in Figure~\ref{F:B} are both rigid braids. Both will give us examples of the transport.  We consider $A$ first.  It has two cycling orbits, $A_1$ and $A_2$, with 3 elements each. Let $s$ be the grey arrow from $A_{1,1}\to A_{2,3}$.  Then the transport of $s$ must be a grey arrow from $\mathbf c(A_{1,1}) = A_{1,2}$ to $\mathbf c(A_{2,3}) = A_{2,1}$, and indeed such an arrow exists and is labeled by the simple element $\sigma_1\sigma_2 = s^{(1)}$.

Example $B$ is more complicated.  Going back to Example~\ref{E:B} we
study the two cycling orbits $B_1 = \{B_{1,1}, B_{1,2}\}$ and $B_{2}
= \{B_{2,1}, B_{2,2}\}$. Since $B_{1,1}$ and $B_{2,1}$ are rigid the
data given in Example~\ref{E:B}  tells us that:}
\setlength{\parskip}{0 pt}
\begin{eqnarray}
\nonumber  B_{1,1} &=& \sigma_2 \sigma_1 \sigma_4 \sigma_3 \sigma_2 \sigma_1 \sigma_5 \sigma_4 \cdot \sigma_2 \sigma_4,  \\
 \nonumber  B_{1,2} &=&  \sigma_2 \sigma_4 \cdot \sigma_2 \sigma_1 \sigma_4 \sigma_3 \sigma_2 \sigma_1 \sigma_5 \sigma_4  \\
  \nonumber  B_{2,2} &=&  \sigma_2 \sigma_4\cdot  \sigma_2 \sigma_1 \sigma_3 \sigma_4 \sigma_3 \sigma_2 \sigma_5 \sigma_4 \\
 \nonumber  B_{2,1} &=&  \sigma_2 \sigma_1 \sigma_3 \sigma_4 \sigma_3 \sigma_2 \sigma_5 \sigma_4 \cdot \sigma_2 \sigma_4.
\end{eqnarray}
{\rm and that the cycling element that takes us from $B_{1,1} \to
B_{2,1}$, that is ${\sigma_2 \sigma_1 \sigma_4 \sigma_3 \sigma_2
\sigma_1 \sigma_5 \sigma_4}$ is the product of the two black arrows
$\sigma_2\sigma_1$ and
$\sigma_4\sigma_3\sigma_2\sigma_1\sigma_5\sigma_4$ in
Figure~\ref{F:B}, both of which are minimal partial cyclings.  Also,
the cycling element that takes us from $B_{2,2} \to B_{2,1}$ is
$\sigma_2\sigma_4$, which is the same as the product of the two
black arrows $\sigma_4$ and $\sigma_2$.   We also see a grey arrow,
call it $s$, from $B_{1,1}\to B_{2,2}$ and its label (which we did
not compute explicitly in Example~\ref{E:B}) is $s = \sigma_1
\sigma_2\sigma_3 \sigma_2 \sigma_1 \sigma_4 \sigma_3 \sigma_2
\sigma_5 \sigma_4 \sigma_3 \sigma_2 \sigma_1.$  The grey arrow from
$B_{1,2}$ to $B_{2,1}$ has the label $t = \sigma_3 \sigma_2 \sigma_1
\sigma_4 \sigma_5 \sigma_4 \sigma_3$. One checks easily that this
part of the diagram is commutative, so that $t = s^{(1)} $ is the
transport of $s$. \setlength{\parskip}{1ex plus 0.5ex minus 0.2ex}

Let's remember that each black arrow represents partial cycling.
This immediately raises a question about the grey arrow that we see
from $B_{3,1}$ to $B_{4,2}$. Is it, too, determined by a commutative
diagram? That leads us to the title of this section: ``Partial
transport".  It will take a while for us to define the concept
precisely. }
\end{example}

We will show that, in the rigid case, there is a
notion of `partial transport' of a minimal simple element,
related to a partial cycling in the same way as the transport is
related to the cycling. More precisely, we will see that a grey arrow
$g$ starting at a vertex $X$ can be {\it partially transported}
along a black arrow $b$, yielding a grey arrow $g^{[b]}$ starting
at $X^b$.  Similarly, we will define the partial transport of $b$
along $g$, which is a black arrow $b^{[g]}$ starting at
$X^{g}$.  We will see that some natural properties are satisfied,
and this will give us more information about the structure of
$\Gamma_X$. In particular, we will show that every two vertices in the same
black component have the same number of grey arrows, and vice versa.

Let $X=\Delta ^p x_1\cdots x_r$ and $Y=\Delta ^p y_1\cdots y_r$ be
two vertices of $\Gamma_X$, where we assume that $X$ is rigid and
$r>1$. Let $b$ be a black arrow going from $X$ to $Y$, that is $X
\stackrel{b}{\longrightarrow} Y$ and $b\preccurlyeq \iota(X)$. We saw
that when we conjugate $X$ by $b$ we are performing a partial
cycling of $X$. Now we want to transport the {\em grey} arrows starting at
$X$ along this partial cycling.

We have $\iota(X)=bc$ for some simple
element $c$. Then $c$ conjugates $Y$ to $\mathbf c(X)$, and these two
elements belong to $USS(X)$. Therefore we can decompose $c=b_2\cdots
b_s$ as a product of minimal simple elements. Denote $b=b_1$. Notice
that $\iota(X)=b_1\cdots b_s$, where $(b_1,\ldots,b_s)$ is an
oriented black path in $\Gamma _X$ by
Corollary~\ref{C:rigid_simple_black_path}.

On the other hand, consider a grey arrow $g $ starting at $X$, and
let $T=X^g $, whose left normal form is $\Delta ^p t_1\cdots
t_r$. We know that its transport $g ^{(1)}$ is defined in such a
way that the following diagram is commutative:
$$
\begin{CD}
X  @>\iota(X)>> \mathbf c(X) \\
@V{g}VV  @VV{g ^{(1)}}V \\
T @>\iota(T)>> \mathbf c(T)
\end{CD}
$$

Moreover, since $g $ is a grey arrow starting at $X$, we have
$\iota(X)\wedge g=1 $ by Proposition~\ref{P:rigid_grey_arrow}.
Hence, by~\cite{gebhardt}, $\iota(X)^{(1)}\wedge g ^{(1)} = 1$.
As $\mathbf c(X)=X^{\iota(X)}$, the definition of the transport
yields $\iota(X)^{(1)}= \iota(\mathbf c(X))$, whence $\iota(\mathbf
c(X))\wedge g^{(1)} =1$. That is, $g ^{(1)}$ is a grey arrow
starting at $\mathbf c(X)$.

Since we have decomposed $\iota(X)$ as a product of black arrows,
the situation is the following:
$$
\begin{CD}
X  @>b_1>>  X_{[1]} @>b_2>>  X_{[2]}  @>b_3>> \cdots  @>b_{s-1}>> X_{[s-1]} @>b_s>>  \mathbf c(X) \\
@Vg VV  @. @. @. @. @VVg ^{(1)}V \\
T  @. @. @. @. @.  \mathbf c(T)
\end{CD}
$$
where $X_{[i]}$ is the conjugate of $X$ by $b_1\cdots b_i$.  Notice that $X_{[1]}=Y$.

We saw in the previous subsection how to transform the black-grey path $(b_s,g ^{(1)})$
into a grey-black path $(g _{s-1}, b'_s)$ such that the following diagram is commutative:
$$
\begin{CD}
X  @>b_1>>  X_{[1]} @>b_2>>  X_{[2]}  @>b_3>> \cdots  @>b_{s-1}>> X_{[s-1]} @>b_s>>  \mathbf c(X) \\
@Vg VV  @. @. @. @VVg_{s-1}V  @VVg ^{(1)}V \\
T  @. @. @. @.  T_{[s-1]} @>b_s'>>  \mathbf c(T)
\end{CD}
$$

We can continue this process, defining for every black-grey path $(b_i,g _i)$ a new
grey-black path $(g _{i-1},b_i')$. At the end we will obtain the following commutative
diagram:
$$
\begin{CD}
X  @>b_1>>  X_{[1]} @>b_2>>  X_{[2]}  @>b_3>> \cdots  @>b_{s-1}>> X_{[s-1]} @>b_s>>  \mathbf c(X) \\
@Vg_0 VV  @VVg_1V        @VVg _2V    @. @VVg_{s-1}V  @VVg ^{(1)}V \\
T  @>b_1'>> T_{[1]} @>b_2'>> T_{[2]}  @>b_3'>> \cdots @>b_{s-1}'>>
T_{[s-1]} @>b_s'>> \mathbf c(T)
\end{CD}
$$

\medskip

\begin{lemma}
In the above situation, $g _0=g $ and $b_1'\cdots
b_s'=\iota(T)$.
\end{lemma}

\begin{proof}
Since $g ^{(1)}$ is a grey arrow starting at $\mathbf c(X)$, we
know by Proposition~\ref{P:rigid_grey_arrow} that $\varphi (\mathbf
c(X))g ^{(1)}$ is a simple element. Since $X$ is rigid, $\varphi
(\mathbf c(X))=\iota(X)$, whence $\iota(X)g^{(1)} $ is simple.
That is, the element $\alpha=b_1\cdots b_sg ^{(1)}$ is simple.

By the above commutative diagram, we know that $\alpha = b_1\cdots
b_s g ^{(1)} = g _0 b_1'\cdots b_s'$, where $g _0$ is a grey
arrow and $(b_1',\ldots,b_s')$ is a black path. Hence, by
Corollary~\ref{C:rigid_simple_grey_black_path}, $g _0 $ and
$b_1'\ldots b_s'$ are determined by $\alpha $. More precisely, the
left normal form of $\varphi(X)\alpha $ is equal to $(\varphi(X)
g_0)(b_1'\cdots b_s')$.

On the other hand, we know by~\cite{gebhardt} that $\alpha =
\iota(X) g ^{(1)}= g\: \iota(T)$. Since $\iota(T)$ is an
element that conjugates $T$ to $\mathbf c(T)$ and since $T$, $\mathbf
c(T)\in USS(X)$, the element $\iota(T)$ can be decomposed as a
product of minimal simple elements $\iota(T)=b_1''\cdots b_t''$.
Moreover, $(b_1'',\ldots,b_t'')$ is an oriented path starting at
$T$, whose associated conjugating element is $\iota(T)$. By
Corollary~\ref{C:rigid_simple_black_path}, $(b_1'',\ldots,b_t'')$ is a
black path. Hence, $(g,b_1'',\ldots,b_t'')$ is a grey-black path, whose associated conjugating
element is
$\alpha $. By Corollary~\ref{C:rigid_simple_grey_black_path}, the
left normal form of $\varphi(X)\alpha $ is equal to $(\varphi(X) g
)(b_1''\cdots b_t'')$. Therefore $g _0=g $ and $b_1'\cdots
b_s'=b_1''\cdots b_t''=\iota(T)$, as we wanted to show.
\end{proof}

Consider the first square of the above commutative diagram:
$$
\begin{CD}
X  @>b>> Y \\
@V{g}VV  @VVg_1 V \\
T @>b'_1>> T_{[1]}
\end{CD}
$$

We want to define $g _1$ as being the {\it partial transport} of
$g $ along $b$. But this will only be consistent if $g _1$ is
completely determined by $X$, $b$ and $g $ (and does not depend
on the decomposition $c=b_2\ldots b_s$ chosen). We will see that
this is the case. Indeed, $g ^{(1)}$ is completely determined by
$X$ and $g $, while $c$ is determined by $X$ and $b$. Hence
$cg^{(1)} $ is determined by $X$, $b$ and $g $. We then use a
decomposition $cg^{(1)} = g _1 b_2'\cdots b_s'$, where $(g
_1,b_2',\ldots,b_s')$ is a grey-black path starting at $Y$. By
Corollary~\ref{C:rigid_simple_grey_black_path}, the grey part of this
path, that is $g _1$, is determined by $cg ^{(1)}$ and
$Y=X^b$, hence it is determined by $X$, $b$ and $g $, as we
wanted to show. More specifically,
Corollary~\ref{C:rigid_simple_grey_black_path} yields $g_1
=\partial(\varphi(Y))\wedge c g ^{(1)}$.

Recalling $bc=\iota(X)$ and $\iota(X) g ^{(1)} =g \:\iota(T)$
\cite{gebhardt}, we can rewrite this as
$$
g _1=\partial(\varphi(Y))\wedge (b^{-1}\iota(X)g ^{(1)})
       =\partial(\varphi(Y))\wedge (b^{-1}g \:\iota(T))
       =\partial(\varphi(X^b)) \wedge (b^{-1} g \:
       \iota(X^{g})).
$$
In particular, this expression for $g_1 $ depends only on $X$, $b$ and
$g$.

\medskip

\begin{definition}
\label{D:partial transport}
{\rm  Let $X\in G$ be rigid and $\ell(X)>1$. Let $b$ and $g$ be arrows in
$\Gamma _X$ starting at $X$, where $b$ is black and $g$ is grey, and let
$Y=X^b$ and $T=X^g$. Then we define the {\it partial transport}
of $g $ along $b$ as follows:}
$$
   g ^{[b]} = \partial(\varphi(Y)) \wedge (b^{-1} g \: \iota(T)).
$$
\end{definition}

Before showing some properties of this partial transport, it will be helpful to also define
the partial transport of a black arrow along a grey one. Using the
analogy between grey and black arrows, considering the graph
$\Gamma_{X^{-1}}$, we just need to recall that
$\partial(\varphi(T^{-1}))=\iota(T)$ and
$\iota(Y^{-1})=\partial(\varphi(Y))$. Hence we have

\begin{definition}
{\rm Let $X\in G$ be rigid and $\ell(X)>1$. Let $b$ and $g$ be arrows in
$\Gamma _X$ starting at $X$, where $b$ is black and $g$ is grey, and let
$Y=X^b$ and $T=X^g$. Then we define the {\it partial transport
of $b $ along $g $} as follows:}
$$
   b^{[g ]} = \iota(T) \wedge (g ^{-1} b \: \partial(\varphi(Y))).
$$
\end{definition}

Fortunately, the partial transport $b^{[g ]}$ is precisely the
element $b_1'$ in the above commutative diagram. This is shown in the
following result:

\begin{lemma}
In the above situation, the following diagram is commutative:
\begin{equation}\label{Diagram:partial_cycling}
\begin{CD}
X  @>b>> Y \\
@V{g}VV  @VVg^{[b]} V \\
T @>b^{[g ]}>> Z
\end{CD}
\end{equation}
\end{lemma}

\begin{proof}
Since $g ^{[b]} = \partial(\varphi(Y)) \wedge (b^{-1} g \:
\iota(T))$, one has $bg ^{[b]}= b\:\partial(\varphi(Y)) \wedge
g \: \iota(T)$. On the other hand, since $b^{[g]}=\iota(T)
\wedge (g ^{-1} b \: \partial(\varphi(Y)))$, one has $g
b^{[g ]}= g \:\iota(T) \wedge  b \: \partial(\varphi(Y))$.
Therefore, $b g ^{[b]} = g b^{[g ]}$, as we wanted to show.
\end{proof}

\noindent {\bf Remark:} Looking at the above diagram, one may be tempted to
conjecture that $bg ^{[b]} = g b^{[g ]} = b\vee g $. This, however, is
{\it not true} in general. One counterexample in the braid group
$B_4$ is given by the rigid braid $X=(\sigma_2)
(\sigma_2\sigma_1\sigma_3\sigma_2) (\sigma_2\sigma_1\sigma_3)
(\sigma_1\sigma_2)$, $g =\sigma _1$ and $b=\sigma _2$. In this
case $b\vee g =\sigma _1\sigma _2\sigma _1$ and $X^{b\vee g}=
(\sigma_1\sigma_3\sigma_2) (\sigma_2\sigma_1\sigma_3\sigma_2)
(\sigma_2) (\sigma_2\sigma_1)$ which is not rigid, that is, does not
belong to $USS(X)$. Since $g b^{[g ]}$ conjugates $X$ to an
element in $USS(X)$, it follows that $b\vee g \neq g b^{[g
]}$. More precisely, in this case $b^{[g ]}= \sigma _2\sigma
_1\sigma _3$, hence $g b^{[g ]} = \sigma _1\sigma _2\sigma
_1\sigma _3 \neq \sigma _1 \sigma _2\sigma _1=b\vee g $.

\bigskip

We will now show some properties of the partial transport which
shall provide more information about the structure of $\Gamma_X$.
First, let us see that the diagram~(\ref{Diagram:partial_cycling})
is the only possible commutative diagram that completes the diagram
$$
\begin{CD}
X  @>b>> Y \\
@V{g}VV   \\
T
\end{CD}
$$

\begin{proposition}\label{P:comm_diagram->partial_cycling}
Let $X\in G$ be rigid and $\ell(X)>1$. Let $b$ be a black arrow in
$\Gamma_X$ starting at $X$ and ending at $Y$ and let $g$ be a
grey arrow starting at $X$ and ending at $T$. Suppose that $b'$ is a
black arrow starting at $T$ and $g'$ is a grey arrow starting at
$Y$, such that $bg'=g b'$. Then $b'=b^{[g]}$ and
$g'=g^{[b]}$.
\end{proposition}

\begin{proof}
By definition, we have  $bg^{[b]} = b\:\partial(\varphi(Y))
\wedge g \:\iota(T)$. However, since $g'$ is a grey arrow
starting at $Y$ one has $g'\preccurlyeq \partial(\varphi(Y))$ so
$bg'\preccurlyeq b\:\partial(\varphi(Y))$. Similarly, since $b'$ is a
black arrow starting at $T$, one has $b'\preccurlyeq \iota(T)$, hence
$bg'= g b' \preccurlyeq g\:\iota(T)$. Therefore, $bg' \preccurlyeq
b g^{[b]}$, that is, $g' \preccurlyeq g^{[b]}$.  Since both
arrows correspond to minimal simple elements, they must be equal,
that is, $g'=g^{[b]}$.

Finally, one has $g b'= b g'= bg^{[b]}=g b^{[g]}$,
hence $b'=b^{[g]}$.
\end{proof}

\begin{corollary}\label{C:black_partial_transport_bijective}
Let $X\in G$ be rigid and $\ell(X)>1$. Let $b$ be a black arrow in
$\Gamma_X$ starting at $X$ and ending at $Y$. Then the partial
transport along $b$ is a bijective map from the set of grey arrows
starting at $X$ to the set of grey arrows starting at $Y$.
\end{corollary}

\begin{proof}
Suppose that $g_1$ and $g_2$ are two grey arrows starting at
$X$ such that $g_1^{[b]}=g_2^{[b]}$. Then $g_1 b^{[g_1]}
= b g_1^{[b]} = b g_2^{[b]} = g_2 b^{[g_2]}$. Hence,
$(g_1,b^{[g_1]})$ and $(g_2,b^{[g_2]})$ are two
grey-black paths starting at $X$ and having the same associated
element. By Corollary~\ref{C:rigid_simple_grey_black_path}, the grey parts of
these paths coincide, so $g_1=g_2$. Hence, the partial
transport along $b$ is injective.

Now let $g'$ be a grey arrow starting at $Y$. Then $(b,g')$ is
an oriented path starting at $X$ formed by a black arrow followed by
a grey arrow. By Proposition~\ref{P:b_rho_->_rho_b}, there exist a
unique grey arrow $g$ and a unique black arrow $b'$ such that
$(g,b')$ is an oriented path starting at $X$ and $bg'=g
b'$. By the above result, one necessarily has $g'=g^{[b]}$.
Hence the partial transport along $b$ is surjective.
\end{proof}

\begin{corollary}\label{C:grey_partial_transport_bijective}
Let $X\in G$ be rigid and $\ell(X)>1$. Let $g$ be a grey arrow in
$\Gamma_X$ starting at $X$ and ending at $T$. Then the partial
transport along $g$ is a bijective map from the set of black
arrows starting at $X$ to the set of black arrows starting at $T$.
\end{corollary}

\begin{proof}
This result is equivalent to the previous one, if we replace $X$ by
$X^{-1}$.
\end{proof}

\begin{corollary}
Let $X\in G$ be rigid and $\ell(X)>1$. Then every two elements in a
black component of $\Gamma_X$ admit the same number of grey arrows
and every two elements in a grey component of $\Gamma_X$ admit the
same number of black arrows.
\end{corollary}

\begin{proof}
If two elements $Y$ and $Z$ belong to the same black component of
$\Gamma_X$, then they can be connected by a sequence of black arrows. By
Corollary~\ref{C:black_partial_transport_bijective}, the partial
transport along each one of these black arrows is a bijection. Hence,
the composition of all these partial transports is a bijection from
the set of grey arrows starting at $Y$ to the set of grey arrows
starting at $Z$.  An analogous proof shows the same statement with reversed
colors.
\end{proof}

In order to complete the picture, we will show that this definition
of partial transport is consistent with the transport defined
in~\cite{gebhardt}. This is done by the following two results.

\begin{lemma}
Let $X$, $Y$, $T$, $b$ and $g $ as above.
$$
\begin{CD}
X  @>b>> Y  @>b^{-1}\iota(X)>> \mathbf c(X) @>b^{(1)}>>  \mathbf c(Y)\\
@V{g}VV  @VVg^{[b]} V @Vg ^{(1)}VV  @VV\left(g^{[b]}\right)^{(1)}V \\
T @>>b^{[g]} > Z @>>(b^{[g]})^{-1}\iota(T)> \mathbf c(T)
@>>\left(b^{[g]}\right)^{(1)}
> \mathbf c(Z)
\end{CD}
$$
\begin{enumerate}
\item The transport of
$g ^{[b]}$ is equal to the partial transport of $g ^{(1)}$
along $b^{(1)}$. That is, $\left(g ^{[b]}\right)^{(1)}=\left(g
^{(1)}\right)^{[b^{(1)}]}$.
\item The transport of $b^{[g]}$ is
equal to the partial transport of $b^{(1)}$ along $g^{(1)}$. That
is, $\left(b ^{[g]}\right)^{(1)}=\left(b
^{(1)}\right)^{[g^{(1)}]}$.
\end{enumerate}
\end{lemma}

\begin{proof}
By~\cite{gebhardt} the transport is multiplicative, whence $(b
g^{[b]})^{(1)} = b^{(1)} (g^{[b]})^{(1)}$. However, since $b
g^{[b]} = g b^{[g]}$, one also has $(b g^{[b]})^{(1)}
=(g b^{[g]})^{(1)} = g^{(1)} (b^{[g]})^{(1)}$. Hence,
$b^{(1)} (g^{[b]})^{(1)} = g^{(1)} (b^{[g]})^{(1)}$. In
other words, if we consider the subdiagrams
$$
\begin{CD}
X  @>b>> Y  @. \hspace{3cm} @. \mathbf c(X) @>b^{(1)}>>  \mathbf c(Y)\\
@V{g}VV  @VVg^{[b]} V  @. @Vg ^{(1)}VV  @VV\left(g^{[b]}\right)^{(1)}V \\
T @>>b^{[g]} > Z  @. @. \mathbf c(T)
@>>\left(b^{[g]}\right)^{(1)}
> \mathbf c(Z)
\end{CD}
$$
then the commutativity of the first one implies the commutativity of
the second one.

The transport is monotonic with respect to $\preccurlyeq$
by~\cite{gebhardt}, whence $b\preccurlyeq\iota(X)$ implies
$b^{(1)}\preccurlyeq\iota(X)^{(1)}=\iota(\mathbf{c}(X))$.  That is, the
transport sends black arrows to black arrows.  As the transport is a
bijection on the set of all arrows and the set of grey arrows is the
complement of the set of black arrows as $X$ is rigid, this means
that the transport sends grey arrows to grey arrows. In particular, $b^{(1)}$ and
$\left(b^{[g]}\right)^{(1)}$ are black arrows and $g^{(1)}$ and
$\left(g^{[b]}\right)^{(1)}$ are grey ones.  By
Proposition~\ref{P:comm_diagram->partial_cycling} and the
commutativity of the second diagram above, this implies $\left(g
^{[b]}\right)^{(1)}=\left(g ^{(1)}\right)^{[b^{(1)}]}$ and
$\left(b ^{[g]}\right)^{(1)}=\left(b
^{(1)}\right)^{[g^{(1)}]}$, as we wanted to show.
\end{proof}

%
%

We will finally show that several partial transports lead to a full transport.  Recall the
commutative diagram
$$
\begin{CD}
X  @>b>>  X_{[1]} @>b_2>>  X_{[2]}  @>b_3>> \cdots  @>b_{s-1}>> X_{[s-1]} @>b_s>>  \mathbf c(X) \\
@Vg VV  @VVg_1V        @VVg _2V    @. @VVg_{s-1}V  @VVg ^{(1)}V \\
T  @>b_1'>> T_{[1]} @>b_2'>> T_{[2]}  @>b_3'>> \cdots @>b_{s-1}'>>
T_{[s-1]} @>b_s'>> \mathbf c(T)
\end{CD}
$$
Recall that $g _1=g ^{[b]}$ and $b_1'=b^{[g ]}$. Denote
$g ^{[b_1,b_2,\ldots, b_i]} = (\cdots ((g ^{[b_1]})^{[b_2]})
\cdots)^{[b_i]}$, where $b_1=b$. Denote also $g _s=g ^{(1)}$.

\begin{lemma}
With the above notations, $g ^{[b_1,\ldots,b_i]}=g _i$ for $i=1,\ldots s$. In
particular, $g ^{[b_1,\ldots,b_s]}=g ^{(1)}$.
\end{lemma}

\begin{proof}
We just need to notice that, by construction, each square in the
above diagram is commutative, and is made of black and grey arrows in
the appropriate sense. Hence by
Proposition~\ref{P:comm_diagram->partial_cycling}, $g_{i}=
g_{i-1}^{[b_i]}$ for $i=2,\ldots,s$, so the result holds.
\end{proof}

Therefore, the transport defined in~\cite{gebhardt} can be seen in the rigid case as
an iterated partial transport along a black path whose associated
conjugating element is $\iota(X)$.  We have thus seen how special
the structure of $\Gamma_X$ is when $X$ is a rigid element with
$\ell(X)>1$.

\section{Periodic elements}
 \label{S:periodic elements}

After the study of rigid elements, we will study other kind of
elements in a Garside group whose associated ultra summit graph is,
in some sense, in the other extreme. This is the case of periodic
elements.

\begin{definition}
{\rm An element $X\in G$ is said to be {\it periodic} if $X^m=\Delta^k$
for some nonzero integers $m$ and $k$.}
\end{definition}

\begin{example} {\rm The $5$-braid $C$ that we gave earlier, in Example~\ref{E:C}, was a periodic braid.    We shall see soon that its graph $\Gamma_C$, depicted in the left sketch in Figure~\ref{F:C-D}, illustrates many of the results in this section. }
\end{example}

In this section we will study the structure of the ultra summit set
of a periodic element $X$. We know that if $X$ were rigid and
$\ell(X)>1$,
then no arrow in $\Gamma_X$ would be bi-colored. The
case of periodic elements is completely the opposite: we will show
that if $X$ is periodic and $\ell(X)>0$, then every arrow in
$\Gamma_X$ is bi-colored. Hence $\Gamma_X = {\cal B}_X= {\cal G}_X$. In
particular, every two elements in $USS(X)$ can be connected by a
sequence of partial cyclings.

\begin{theorem} \label{T:X periodic imples ell(X)=1} {\rm (\cite{bestvina} for Artin-Tits groups, \cite{CMW} for Garside groups.)}
If $X\in SSS(X)$ is periodic and it is not a power of $\Delta$, then
$\ell(X)=1$.
\end{theorem}

\begin{remark}   {\rm A new proof of Theorem~\ref{T:X periodic imples ell(X)=1} for all Garside groups  follows immediately from the decomposition that is given in Theorem 2.9 of \cite{BGGM-I}. }
\end{remark}
\bigskip

Recall that an {\it atom} is an indivisible simple element, and that
every Garside group $G$ is generated by its atoms. We will show
that if an element $X\in USS(X)$ is periodic, then all its minimal
simple elements are atoms and correspond to partial cyclings. We
will also determine which atoms are minimal simple elements for a
given periodic element.

Here is the main result concerning periodic elements.

\begin{proposition}\label{P:mse_for_periodic}
Let $X\in USS(X)$ be a periodic element which is not a power of
$\Delta$. The set of minimal simple elements for $X$ is equal to the
set of atoms dividing $\iota(X)$ and $\iota(X^{-1})$, that is, the
left descent set of $\iota(X)\wedge \iota(X^{-1})$, in the
terminology of Artin-Tits groups.
\end{proposition}

\begin{proof}
It is well known~\cite{el-m} that $\ell(Y)=\ell(Y^{-1})$ for every
$Y\in G$. This implies that $SSS(Y^{-1})$ consists of the inverses
of the elements in $SSS(Y)$. In our case, since $X$ is periodic and
it is not a power of $\Delta$, we have $\ell(X)=1$. Hence, $USS(X)=SSS(X)$
and $USS(X^{-1})=SSS(X^{-1})=\{ Y^{-1}\mid Y\in SSS(X)\}$ by
Lemma~\ref{L:periodic_length_one}.
Moreover, a simple element $t$
satisfies $X^t\in USS(X)$, that is $\ell(t^{-1}Xt)=1$, if and only
if $\ell(t^{-1}X^{-1}t)=1$, which is equivalent to
$(X^{-1})^t \in USS(X^{-1})$. Therefore, the minimal simple elements
for $X$ and for $X^{-1}$ coincide.

Let $\Delta^p x_1$ be the left normal form of $X$, and recall that
$\iota(X)=\tau^{-p}(x_1)$ and $\iota(X^{-1})=\partial(x_1)$. Let $s$
be a minimal simple element for $X$ (thus for $X^{-1}$). By
Corollary~\ref {C:prefix}, we have $s\preccurlyeq \iota(X)$
or $s\preccurlyeq \iota(X^{-1})$. Replacing $X$ by $X^{-1}$ if necessary,
we can assume that $s\preccurlyeq \iota(X^{-1})=\partial(x_1)$. This
implies that $x_1 s$ is simple.

We will show that we also have $s\preccurlyeq \iota (X) = \tau^{-p}(x_1)$
or, equivalently, that $\tau^p(s)\preccurlyeq x_1$. We know that
$\tau^p(s)\vee x_1 = x_1t$, for some simple element $t$. We must
then show that $t=1$.

Notice that $\ell(s^{-1}Xs)=\ell(s^{-1}\Delta^p x_1 s) = 1$, hence
$\:\tau^p(s)^{-1}x_1s\:$ is simple. That is, $\tau^p(s)\preccurlyeq
x_1s$. Then, $x_1s$ is a multiple of $\tau^p(s)$ and of $x_1$, so it
follows that $x_1 t \preccurlyeq x_1 s$, thus $t\preccurlyeq s$. This yields
$\tau^p(t)\preccurlyeq \tau^p(s)\preccurlyeq x_1t$, whence $\tau^p(t)^{-1}x_1
t $ is simple, implying $\ell (t^{-1}Xt)=1$.  By minimality of $s$,
it follows that either $t=1$ (as we want to show) or $t=s$.

Suppose that $t=s$, that is, $\tau^p(s)\vee x_1 = x_1s$. We will see
that in this case $X$ cannot be periodic, since we will show by
induction that $\ell(X^k)>0$ and that $s\not\preccurlyeq \iota(X^k)$, for
every $k>0$. For $k=1$, we know that $\ell(X)=1$, and if we had
$s\preccurlyeq \iota(X) =\tau^{-p}(x_1)$, then $\tau^p(s)\preccurlyeq x_1$, so
$\tau^p (s)\vee x_1 = x_1 = x_1s$, which would imply that $s=1$ in
contradiction to our assumption. Hence, the claim is true for $k=1$.

Suppose the claim is true for some $k\geq 1$. Let $\Delta^q
z_1\cdots z_m$ be the left normal form of $X^k$, where $m\geq 1$.
Then $X^{k+1}=\Delta^p x_1 \Delta^q z_1 \cdots z_m $. In order to
show that $\ell(X^{k+1})>0$ it suffices to prove that
$x_1\tau^{-q}(z_1)$
is not a multiple of $\Delta$. This follows if we
can show that $\partial(x_1) \not \preccurlyeq \tau^{-q} (z_1)$, that is,
that $\iota(X^{-1})\not \preccurlyeq \iota(X^k)$. The latter is true
since we are assuming that $s\preccurlyeq \iota(X^{-1})$, while by
induction hypothesis $s\not \preccurlyeq \iota(X^k)$. Hence
$\ell(X^{k+1})>0$. Moreover, in this case $\iota(X^{k+1})$ is the
maximal simple prefix of $\tau^{-p}(x_1)\tau^{-p-q}(z_1)$.  If we
had $s\preccurlyeq \tau^{-p}(x_1)\tau^{-p-q}(z_1)$, then
$\tau^p(s)\preccurlyeq x_1 \tau^{-q}(z_1) = x_1 \:\iota(X^k)$. In this
case, $x_1 \:\iota(X^k)$ would be a multiple of $\tau^p(s)$ and of
$x_1$, whence we would have $x_1 s \preccurlyeq x_1 \iota(X^k)$ as we
assume $\tau^p(s)\vee x_1 = x_1s$. This is not possible, since by
induction hypothesis $s\not\preccurlyeq \iota(X^k)$. Therefore,
$s\not\preccurlyeq\tau^{-p}(x_1)\tau^{-p-q}(z_1)$, whence $s\not\preccurlyeq
\iota(X^{k+1})$ which proves the claim.

We have then shown that if $t=s$, the element $X$ would not be
periodic. Hence, we have $t=1$. This means that $\tau^p(s)\vee x_1
=x_1$, that is, $\tau^p(s)\preccurlyeq x_1$, which implies that $s\preccurlyeq
\tau^{-p}(x_1)=\iota(X)$.  Therefore, every minimal simple element
$s$ for $X$ must be a prefix of $\iota(X^{-1})$ and also of
$\iota(X)$.

Now choose any atom $a\preccurlyeq s$. We know that $a\preccurlyeq
\iota(X)=\tau^{-p}(x_1)$ and that $x_1a$ is simple. Therefore,
$a^{-1}Xa$ has canonical length 1. By minimality of $s$, we have
$a=s$. Hence, every minimal simple element for $X$ is an atom
dividing both $\iota(X)$ and $\iota(X^{-1})$.

Conversely, an atom $a$ which divides both $\iota(X)$ and
$\iota(X^{-1})$ satisfies $\tau^p(a)\preccurlyeq x_1$ and $x_1a$ is
simple, whence $\ell(a^{-1}Xa)=1$. Since an atom has no proper
prefixes, it is a minimal simple element for $X$.
\end{proof}

The most important consequence of this result is that one can
connect every two elements in the ultra summit set of a periodic
element by a sequence of partial cyclings.

\begin{corollary}\label{C:periodic_main}
Let $X\in USS(X)$ and $Y\in USS(Y)$ be periodic elements. Then
$X$ and  $Y$ are conjugate if and only if there exists a sequence of
elements $X=X_1, X_2, \ldots, X_m=Y$, such that for all
$i=1,\ldots,m-1$, the element $X_{i+1}$ is a partial cycling of $X_i$
and $X_i\in USS(X_i) = USS(X)$.
\end{corollary}

\begin{proof} If $X=\Delta^m$ for some $m$, then $USS(X)=\{X\}$, so the
result is trivially true. We can then assume that $X$ is not a power
of $\Delta$.

Suppose that $X$ and $Y$ are conjugate.
By~\cite{gebhardt} there exists a chain of elements $X=X_1, X_2,
\ldots, X_m=Y$, such that $X_{i+1}=X_i^{s_i}$, where $X_i\in
USS(X)$ and $s_i$ is a minimal simple element for $X_i$, for
$i=1,\ldots,m-1$.
By Proposition~\ref{P:mse_for_periodic}, each $s_i$ is an atom
dividing $\iota(X_i)$, whence $X_{i+1}$ is a partial cycling of
$X_i$, for $i=1,\ldots,m-1$.

The converse is trivial.
\end{proof}

\bigskip

We will see in Proposition~\ref{P:going to USS by partial} that, given $X\in G$, we can obtain some
$X'\in USS(X)$ by a finite number of partial cyclings. Given
$X,Y\in G$ which are conjugate, we can hence obtain $X'\in USS(X)$ and
$Y'\in USS(X)$ using partial cyclings. If $X$ is periodic, we have
also shown that one can go from $X'$ to $Y'$ by a finite number of
partial cyclings.

But in general, even if $X$ is periodic, one cannot go from $X$ to
$Y$ by partial cyclings, since a partial cycling will never decrease
the infimum of an element, and $Y$ does not necessarily have
maximal infimum.

\section{Applications}
\label{S:applications}

\subsection{Complexity of the CDP/CSP for pseudo-Anosov and periodic braids}
\label{SS:complexity of the CDP/CSP for PA and periodic braids}
The work in this paper is applicable to all Garside groups, and no assumptions are made that restrict attention to the braid groups $B_n$.   In this section we show, by applying our results to braids, that we have made progress toward one of our long-range goals.  That goal is to find a solution to the CDP/CSP for PA elements in $B_n$ which is {\it polynomial}, that is  bounded above by a polynomial  in both $n$ and $\ell(X)$.  In this section we discuss our progress, at this time, toward that goal, and what remains to be done.

Putting $X$ into left normal form is known (see \cite{el-m, BKL2001}) to have polynomial complexity.  It follows from the main result in \cite{BKL2001} that there is a polynomial bound to how many times one must cycle and decycle an element $X$ in order to bring $X$ to a representative which is in $SSS(X)$, so without loss of generality we may assume that $X$ is in left normal form and in $SSS(X)$. But at this time we do not know how many times one must cycle to bring $X$ into $USS(X)$.  That problem may yield, using the machinery developed in $\S$2 of \cite{BGGM-I}, however at this writing it remains an open question.  Let us assume from now on that $X\in USS(X)$.

\large
{\bf  (A)  The PA case:}
\normalsize
We have proved in \cite{BGGM-I} that, if $X\in USS(X)$ is PA, then there is an integer $m$ which is at most $((n)(n-1)/2)^3$, such that every element in $USS(X^m)$ is rigid if we use the usual Garside structure.  If we use the BKL structure \cite{BKL1998}  then $m$ is bounded by $(n-1)^3$.  In view of the fact, proved in \cite{gm}, that roots of PA braids are unique, this means that, in searching for a polynomial solution to the CDP/CSP for PA braids, we may assume that they are rigid.  Since any rigid element $X$  of a Garside group is in $USS(X)$, this means we may assume that $X$ is PA and that $X\in USS(X)$ is rigid.

We have computed many many examples of $USS(X)$ when $X$ is PA and rigid, using random
searches, and on the basis of the evidence found that in the generic
case $USS(X)$ has either 2 orbits, where one is the conjugate of the
other by $\Delta$, or 1 orbit which is conjugate to itself by
$\Delta$.   Those are the two cases that we saw in Example~\ref{E:A} and Figure~\ref{F:A}. We found this behavior over and over again, in calculations with very large numbers of randomly chosen examples, when we restricted our attention to PA braids that are rigid.  A question that remains, for future work, is: If $X$ is a PA and rigid element in $B_n$, is the size of $USS(X)$ bounded above by some polynomial in $n$ and $\ell(X)$?  This is where our structure theorems about $USS(X)$ become very important.  We present two examples which illustrate the problems that remain to be solved in the PA, rigid case.  Let's remember that in such cases we understand cycling very well, but the partial cyclings and partial twisted decyclings that connect cycling orbits present new combinatorial challenges.

We would like to give some examples to illustrate the difficulties, but  encounter a problem on how to both present interesting examples and draw good pictures, when $USS(X)$ is unexpectedly large.  To overcome the difficulty,  recall that for every $X\in G$ there is an automorphism $\phi$ of $\Gamma_X$ that sends every vertex $Y$ to $\mathbf c(Y)$ and every arrow $s$ to its transport $s^{(1)}$.  We can hence define a quotient graph $\Gamma_X/ \phi$ whose vertices are the cycling orbits in
the vertex set of $\Gamma_X$. Recall also that $\phi$ preserves the colors of the arrows, hence the arrows in $\Gamma_X/\phi$ also have a well-defined color (black or grey).  We can then draw the quotient graphs $\Gamma_X/ \phi$ to have an idea how the distinct orbits are connected in the graph $\Gamma_X$.  In the two quotient graphs, Figures~\ref{F:E} and \ref{F:F} below, the vertex label $E_i$ or $F_i$ means a cycling orbit.

\begin{example} \label{E;E}
{\rm Figure~\ref{F:E} illustrates one of the difficulties.   We consider the 12-braid:
\begin{eqnarray}
\nonumber E = &(\sigma_2 \sigma_1 \sigma_7 \sigma_6 \sigma_5 \sigma_4 \sigma_3 \sigma_8 \sigma_7 \sigma_{11} \sigma_{10}) \cdot (\sigma_1 \sigma_2 \sigma_3 \sigma_2 \sigma_1 \sigma_4 \sigma_3 \sigma_{10}) \cdot  \\
\nonumber & (\sigma_1 \sigma_3 \sigma_4 \sigma_{10}) \cdot  (\sigma_1 \sigma_{10}) \cdot  (\sigma_1 \sigma_{10} \sigma_9 \sigma_8 \sigma_7 \sigma_{11}) \cdot  (\sigma_1 \sigma_2 \sigma_7 \sigma_{11})
\end{eqnarray}
The braid $E$ is a PA, rigid braid, with cycling orbits of length 6. It turns out that  $\Gamma_E$ has 264 elements, but the quotient graph $\Gamma_E/ \phi$  has 44 vertices.  Of course we cannot distinguish between elements in the same cycling orbit in the quotient graph, however with some extra information we will be able to understand $\Gamma_E$ too:} \end{example}
\be
\item  There is a black arrow from $E_{i,j}$ to $E_{i+2,j}$ for every $i\neq 43,44$.
\item There is a black arrow from $E_{i,j}$ to $E_{i-2,j+1}$ for every $i\neq 1,2$.
\item The product of two consecutive black arrows in opposite senses is trivial on orbits, but (in this example) corresponds to cycling of the elements in the orbit.
\item There is a grey arrow from $E_{i,j}$ to $E_{i+1,j-1}$ if $i$ is odd.
\item  There is a grey arrow from $E_{i,j}$ to $E_{i-1,j-1}$ if $i$ is even.
\item  Every grey arrow corresponds to a twisted decycling. (Hence $\tau(E_i)=E_{i+1}$ for every odd  $i$).
\ee
\begin{figure}[htpb]
\centerline{\includegraphics[scale=1.0]{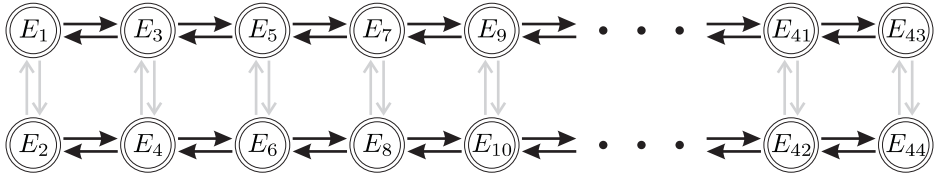}}
\caption{The quotient graph $\Gamma_E/ \phi$ illustrates a very long path of iterated partial cyclings.}
\label{F:E}
\end{figure}

\begin{example} \label{E;F}
{\rm Figure~\ref{F:F} illustrates another way in which partial cycling leads to difficult combinatorial problems. Consider the 12-braid $F$:
\begin{eqnarray}
\nonumber F = & (\sigma_3 \sigma_2 \sigma_1 \sigma_4 \sigma_6 \sigma_8 \sigma_7 \sigma_6 \sigma_9 \sigma_{10}\sigma_{11} \sigma_{10}) \cdot (\sigma_1 \sigma_2 \sigma_4 \sigma_3 \sigma_2 \sigma_1 \sigma_5 \sigma_7 \sigma_{10} \sigma_{11} \sigma_{10}) \cdot  \\
\nonumber  &  (\sigma_3 \sigma_5 \sigma_7
\sigma_{10 }\sigma_{11} \sigma_{10}) \cdot (\sigma_3 \sigma_5 \sigma_7 \sigma_6 \sigma_8 \sigma_{10} \sigma_{11})
\end{eqnarray}

\begin{figure}[htpb]
\centerline{\includegraphics[scale=1.0]{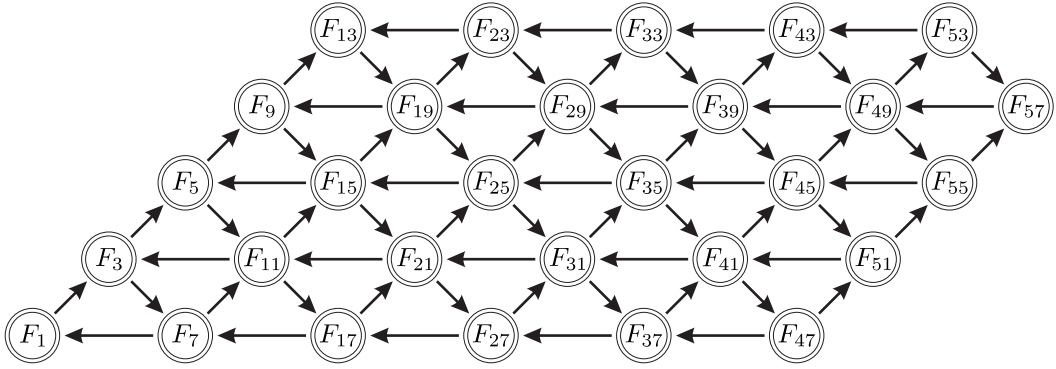}}
\caption{The quotient graph $\Gamma_F/ \phi$ illustrates multiple minimal partial cyclings at a vertex.}
\label{F:F}
\end{figure}

Since $F$ is rigid and since $\ell(F) = 4$ each cycling orbit has 4 elements. The graph $\Gamma_F$ has two isomorphic black components, each made of 29 cycling orbits.  We show one of them.  Its graph  has the following properties, the key one being Property 3:}
\end{example}
\be
\item Every black arrow pointing toward the right (joining orbits $i$ and $j$) goes from $F_{i,k}$ to $F_{j,k}$.
\item Every black arrow pointing to the left (joining orbits $i$ and $j$) goes from $F_{i,k}$ to $F_{j,k+1}$.
\item The concatenation of three black arrows forming a `small triangle' corresponds to a cycling. Hence, by Lemma~\ref {L:decomposition on partial cycling},  the initial factor of every vertex of $\Gamma_F$ which is not in $F_1$ can be decomposed in several ways as a product of minimal simple elements.
For example, the element $F_{11,1}$  has 6 different
decompositions of $\iota(F_{11,1})$, given by the small triangles 11-3-7,
11-3-5, 11-15-5, 11-15-21, 11-17-21, 11-17-7.
This causes `branching' at the vertex $F_{11,1}$.
\item We have grey arrows, corresponding to partial twisted decyclings, going from $F_{2k-1,j}$ to $F_{2k,j-1}$ and also from $F_{2k,j}$ to $F_{2k-1,j-1}$.
\ee

Here are the problems that remain, with regard to PA braids. First, we need to learn how many times we must cycle to bring $X\in SSS(X)$ into a closed cycling orbit, i.e. into $USS(X)$. This is `open problem 3' of $\S$1.4 of \cite{BGGM-I}.   We have shown (in \cite{BGGM-I}) that once we are in $USS(X)$, we may assume that $X$ is rigid.  If $X$ is in left normal form and is rigid, then the length of every cycling orbit is either $\ell(X)$ or $2\ell(X)$. Lemma~\ref{L:decomposition on partial cycling} shows that very very special combinatorial conditions are required for partial cycling, and even more for iterated partial cycling,  to occur.  Yet Examples $E$ and $F$ show that $USS(X)$ can have surprises.  We must  find polynomial bounds on (i) the lengths of paths in $\Gamma_X$ and (ii) the number of such paths (that is, the combinatorics introduced by branching).  For fixed braid index and length $\ell(X)$ we need a universal bound on the lengths and numbers of such paths.  If we can solve all these problems, then we should be able to solve Open Question 2 of $\S$1.4 of \cite{BGGM-I}, that is to find a polynomial bound on $|USS(X)|$ and so (using all our other work in \cite{BGGM-I} and this paper, as well as drawing heavily on the literature), obtain a polynomial solution to the CDP/CSP for PA braids.

\large
{\bf (B) The periodic case:}
\normalsize
The situation for periodic braids is quite different from that for PA braids.  On the one hand, periodic braids are quite simple. It is well-known, from the work of Eilenberg \cite{Eilenberg} and
Ker\'ekj\'art\'o \cite{Ker} that every periodic element  $X\in B_n$ is conjugate to either a power of $\delta = \sigma_{n-1}\sigma_{n-2}\cdots\sigma_1$ or a power of $ \varepsilon = \delta\sigma_1$.  However, unfortunately, it turns out that $|USS(\delta)| = 2^{n-2}$ and $|USS(\varepsilon)| = (n-2)2^{n-3}$ (this will be proved in \cite{BGGM-III}), so without even considering powers of $\delta$ and $\varepsilon$ we have exponential growth.
We will nevertheless arrive at a polynomial solution to the CDP/CSP in \cite{BGGM-III}, by introducing new tricks that  put the two known Garside structures on $B_n$ to work,  and make use of Garside structures on other Artin groups too.   That is the content of the manuscript \cite{BGGM-III}.

\subsection{A new solution to the CDP/CSP in Garside groups}
\label{SS:a new solution to CDP/CSP in Garside groups}

In this section, we use earlier results to present a new
solution of the CDP/CSP problem in a Garside group $G$, that is,
we give an algorithm which determines whether two elements $X,Y\in G$ are
conjugate. Moreover, if $X$ and $Y$ are conjugate, the algorithm
finds a conjugating element $\alpha$ such that $X^{\alpha}=Y$.

Recall that the algorithm in~\cite{gebhardt} is based on computing
$USS(X)$ and one element $Y'\in USS(Y)$. Then, $X$ and $Y$ are
conjugate if and only if $Y'\in USS(X)$. Hence, if one uses the
algorithm in~\cite{gebhardt}, one must compute the entire set $USS(X)$.
Corollary~\ref{C:components_intersect}, however, shows that this can
be avoided.  Here is the promised algorithm.  It is based upon Corollary~\ref{C:new criterion}

\bigskip
\noindent {\sc Algorithm 3.}

\medskip
 {\bf Input:} $X,Y\in G$.

\medskip {\bf Output:} $\alpha\in G$ such that $X^{\alpha}=Y$, or `Fail' if
$X$ and $Y$ are not conjugate.

\begin{enumerate}

 \item Using cyclings and decyclings, compute $X'\in USS(X)$ and
 $a\in G$ such that $X^a=X'$.

 \item Using cyclings and decyclings, compute $Y'\in USS(Y)$ and
 $b\in G$ such that $Y^b=Y'$.

 \item Using Algorithm 1, compute ${\cal B}_{X'}$ and, for each vertex
  $V$ of ${\cal B}_{X'}$, an element $c_{(X',V)}$ conjugating $X'$ to $V$.

 \item Using Algorithm 2, compute ${\cal G}_{Y'}$ and, for each vertex
  $V$ of ${\cal G}_{Y'}$, an element $c_{(Y',V)}$ conjugating $Y'$ to $V$.

 \item If ${\cal B}_{X'}\cap {\cal G}_{Y'}=\emptyset$ return `Fail'.

 \item Choose $V\in {\cal B}_{X'}\cap {\cal G}_{Y'}$ and return
 $\;a\: c_{(X',V)} \:c_{(Y',V)}^{-1}\: b^{-1}$.

\end{enumerate}

\bigskip

Therefore, in order to determine whether $X$ and $Y$ are conjugate,
we just need to compute one black component and one grey component.
The union of these two sets is in general smaller than the whole
ultra summit set.  Moreover, this procedure provides a conjugating
element, since the algorithm computes graphs in such a way that we
know how to join any two vertices in our graphs, and each path in
the graph yields a conjugating element between the initial and the
final vertex.

It should be mentioned that, in the worst case, this new algorithm is not
better than the one in~\cite{gebhardt}, since there are examples in
which either ${\cal B}_{X'}=USS(X)$ or ${\cal G}_{Y'}=USS(Y)$. This
is the case, for instance, for the periodic elements treated in
Section~4.

\subsection{Partial cycling subsumes decycling}
\label{SS:partial cycling subsumes decyling}

In Corollary~\ref{C:partial cycling and USS(X)} we showed that we only need partial cyclings and partial
twisted decyclings to connect a pair of elements in the same
ultra summit set.  We now prove that, given any $X\in G$, where in general $X\notin USS(X)$, we can conjugate
$X$ to an element in $USS(X)$ by applying a finite number of partial
cyclings.

We know by~\cite{el-m} (see also the review in \cite{BGGM-I}) that we can obtain an element in $USS(X)$ from any $X\in G$ by iterated cyclings and decyclings. It is
clear that a cycling is a particular case of a partial cycling, but
we will now prove that a decycling can also be seen, in some sense,
as a composition of several partial cyclings, provided that the
element involved has maximal infimum in its conjugacy class.

\begin{lemma}\label{L:decycling_is_partial}
Let $X\in G$, and let $p=\inf(X)$ and $r=\ell(X)$. If the infimum of
$X$ is maximal in its conjugacy class, then there exists a sequence
$X = X_1, X_2, \ldots, X_r = \mathbf d(\tau^{-p}(X))$, such that
$X_{i+1}$ is a partial cycling of $X_i$, for $i=1,\ldots,r-1$.
\end{lemma}

\begin{proof}
For simplicity, we will consider $Y=\tau^{-p}(X)$, and we will show
that there is a sequence $\tau^p(Y)= Y_1, Y_2,\cdots, Y_r = \mathbf
d(Y)$, consisting of partial cyclings. The claim is trivial for
$r=0$, so assume $r>0$ and let $\Delta^p y_1\cdots y_r$ be the left
normal form of $Y$. By definition, $\mathbf d(Y)=y_r \Delta^p y_1
\cdots y_{r-1}=Y^{y_r^{-1}}$. Since $Y$ commutes with itself, one
also has $\mathbf d(Y)= Y^{(Y y_r^{-1})} = Y^{(\Delta^p y_1\cdots
y_{r-1})} = (\tau^p(Y))^{(y_1\cdots y_{r-1})}$.

Let $Y_1=\tau^p(Y)$ and $Y_{i+1}=Y_i^{y_i}$, for $i=1,\ldots,r-1$
(thus $Y_r=\mathbf d(Y)$). We will show that $Y_{i+1}$ is a partial
cycling of $Y_i$ for $i=1,\ldots,r-1$.

Notice that for every element $Z$ such that $\inf(Z)=p$,
the first factor in the left normal form of
$Z\Delta^{-p}$ (whose infimum equals $0$) is $\iota(Z)$.
 In other words, a simple
element $s$ performs a partial cycling on $Z$ if and only if $s$ is
a prefix of $Z \Delta^{-p}$.

Notice also that $Y_1 \Delta^{-p} = \tau^p(Y) \Delta^{-p} =
\Delta^{-p} Y = y_1\cdots y_r$ and
that one has $Y_i \Delta^{-p} =
Y_1^{(y_1\cdots y_{i-1})}\Delta^{-p} = (y_1\cdots y_{i-1})^{-1} (Y_1
\Delta^{-p}) \tau^{-p}(y_1\cdots y_{i-1}) = y_i \cdots y_r
\tau^{-p}(y_1\cdots y_{i-1})$ for every $i=2,\ldots,r-1$.  Moreover,
$\inf(Y_i)=p$, since it cannot be greater by hypothesis, and it
cannot be smaller since $Y_i\Delta^{-p}$ is positive.

This implies that for $i=1,\ldots,r-1$ the simple factor $y_i$ is a
prefix of $Y_i \Delta^{-p}$, where $\inf(Y_i)=p$. Therefore,
$Y_{i+1}=Y_i^{y_i}$ is a partial cycling of $Y_i$ for
$i=1,\ldots,r-1$.
\end{proof}

\begin{proposition}\label{P:going to USS by partial}
Given $X\in G$, there exists a sequence $X=X_1,X_2,\ldots, X_k=Y$ in
$G$, such that $Y\in USS(X)$, and $X_{i+1}$ is a partial cycling of
$X_i$ for $i=1,\ldots, k-1$.
\end{proposition}

\begin{proof}
It is well known~\cite{el-m} that, by applying a finite sequence of
cyclings to an element in $G$, one can achieve the maximal infimum
in its conjugacy class. Since cyclings are a particular case of
partial cyclings, we can assume that $X$ has maximal infimum in its
conjugacy class.

By~\cite{el-m} and \cite{gebhardt}, we know that $\mathbf c^s
(\mathbf d^t(X))\in USS(X)$, for some positive integers $s$ and $t$.
On the other hand, we know that the ultra summit set is invariant
under $\tau$. Hence, if we denote $p=\inf(X)$, we have
$Y=\tau^{-pt}(\mathbf c^s (\mathbf d^t(X)))\in USS(X)$.

By Lemma~\ref{L:tau_c_d}, $\tau$ commutes with cycling and
decycling, hence we can write $Y=\mathbf c^s ((\mathbf d\circ
\tau^{-p})^t(X))\in USS(X)$. Finally, by
Lemma~\ref{L:decycling_is_partial}, together with the fact that $X$ has maximal
infimum in its conjugacy class, one can decompose every $(\mathbf d
\circ \tau^{-p})$ as a product of partial cyclings. Since the final
application of cyclings also correspond to partial cyclings, the
result follows.
\end{proof}

\subsection{An application to the theory of reductive groups}
\label{SS:an application to the theory of reductive groups}

In $\S$\ref{S:periodic elements} we showed that if $X$ is periodic and $\ell(X)>0$, then every arrow in
$\Gamma_X$ is bi-colored. Hence $\Gamma_X={\cal B}_X= {\cal G}_X$. In
particular, every two elements in $USS(X)$ can be connected by a
sequence of partial cyclings.

It was communicated to us by Jean Michel that this has important
consequences in the theory of reductive groups. Namely, in a
reductive group, the Deligne-Lusztig varieties are related to the
elements of the braid monoid of the Weil group $W$. In particular,
the Deligne-Lusztig varieties which appear in Brou\'e's conjecture
correspond to periodic elements in an Artin-Tits group. There is a
conjecture which states that there is essentially a unique
Deligne-Lusztig variety for each period.   Fran\c{c}ois
Digne  and Jean Michel have shown that, if two periodic elements are conjugate by a
partial cycling, then their corresponding Deligne-Lusztig varieties
are essentially the same. In the case of Artin braid groups (which
corresponds to the linear algebraic groups), it is known that two
periodic elements of the same period are conjugate \cite{gm}. Hence, the
results of Section~\ref{S:periodic elements} show that the above conjecture is true, at
least for Artin braid groups (Artin-Tits groups of type A). In order
to show the conjecture in general, it remains to be shown that every
two periodic elements of the same period, in any spherical type
Artin-Tits group, are conjugate. We refer to \cite{DM} for details on this problem.

\footnotesize
\begin{tabular}{lll}
 {\bf Joan S. Birman}        &  {\bf Volker Gebhardt}  &  {\bf Juan Gonz\'alez-Meneses} \\
 Department of Mathematics, & School of Computing and Mathematics, & Departamento de \'Algebra, \\
 Barnard College andColumbia University, & University of Western Sydney, & Universidad de Sevilla,\\
 2990 Broadway, & Locked Bag 1797, &  Apdo. 1160, \\
 New York, New York 10027, USA. & Penrith South DC NSW 1797, Australia, & 41080 Sevilla, Spain.\\
 {\tt jb@math.columbia.edu} & {\tt v.gebhardt@uws.edu.au} & {\tt meneses@us.es}
\end{tabular}


\begin{thebibliography}{99}

\bibitem{bestvina} V. Bestvina. {\it Non-positively curved aspects of Artin groups of finite
type.} Geometry and Topology, 3 (1999) 269-302.

\bibitem{BGGM-I}  J. Birman, V. Gebhardt and J. Gonz\'alez-Meneses, {\it Conjugacy in Garside groups I:
Cycling, Powers and Rigidity}, preprint arXiv math.GT/0605230.

\bibitem{BGGM-III}   J. Birman, V. Gebhardt and J. Gonz\'alez-Meneses, {\it Conjugacy in Garside groups III:Periodic Braids}, manuscript in preparation.

\bibitem{BKL1998} J. Birman, K.Y.Ko and S.J.Lee, {\it A new approach to the word and conjugacy problems in the braid groups}, Adv. Math. {\bf 139}, No. 2, (1998), 322-353.

\bibitem{BKL2001}  J. Birman, K.Y. Ko and S.J. Lee,   {\it The infimum, supremum and geodesic length of a braid conjugacy class}, Adv. Math. (2001), {\bf 164}, No. 1, (2001), 41-56.

\bibitem{B-S} E. Brieskorn and K. Saito, {\it Artin-Gruppen und Coxeter-Gruppen}, Invent. Math. {\bf 17} (1972), 245-272.

\bibitem{CMW} R. Charney, J. Meier and K. Whittlesey. {\it Bestvina's normal form complex
and the homology of Garside groups.} Geometriae Dedicata 105 (2004),
171-188.

\bibitem{D-P}  P. Dehornoy  and L. Paris, {\it Gaussian groups and Garside groups, two generalizations of Artin
groups}, Proc. London Math. Soc. {\bf 79} (1999), No. 3, 569-604.

\bibitem{DM} F. Digne and J. Michel, {\it Endomorphisms of Deligne-Lusztig varieties}, preprint  arXiv math.RT/0509011

\bibitem {Eilenberg} S. Eilenberg, {\it Sur les transformations p$\acute{e}$riodiques de la surface de la sph$\acute{e}$re}, Fund. Math. {\bf 22} (1934), 28-41.

\bibitem{el-m} E. ElRifai and H. Morton, {\it Algorithms for positive braids}, Quart. J. Math. Oxford
Ser (2), {\bf 45} (180) (1994), 479-497.

\bibitem{F-GM} N. Franco and J. Gonz\'alez-Meneses, {\it The conjugacy problem for braid groups and Garside
groups}, Journal of Algebra, {\bf 266}, No. 1 (2003), 112-132.

\bibitem{garside} F. Garside, {\it The braid group and other groups}, Quart. J. Math Oxford {\bf 20}
(1969), 235-254.

\bibitem{gebhardt} V. Gebhardt, {\it A new approach to the conjugacy problem in Garside
groups}, Journal of Algebra, {\bf 292}, No. 1 (2005), 282-302.

\bibitem{gm} J. Gonz\'alez-Meneses, {\it The $n^{th}$ root of a braid is unique up to conjugacy}, Algebraic
and Geometric Topology {\bf 3} (2003), 1103-1118.

\bibitem{GM-W} J. Gonz\'alez-Meneses and B. Wiest, {\it On the structure of the centralizer of a braid}, Ann. Sci. \'Ecole Norm. Sup. (4) {\bf 37} (2004), no. 5, 729-757.

\bibitem{Ker}  B. de Ker\'ekj\'art\'o,  {\it $\ddot{U}$ber die periodischen Transformationen der Kreisscheibe und der Kugelfl$\ddot{a}$che}, Math. Annalen {\bf 80} (1919), 3-7.

\bibitem{Michel} J. Michel, {\it A note on words in braid monoids}, J. Algebra {\bf 215} (1999), no. 1, 366-377.

\bibitem{thurston} W. Thurston, {\it On the topology and geometry of diffeomorphisms of
surfaces}, Bull. Amer. Math. Soc. {\bf 19} (1988), 109-140.

\end{thebibliography}
\end{document}